\documentclass[11pt,reqno]{amsart}
\usepackage{amsmath,amsthm,amssymb,verbatim,enumerate,ifthen}
\usepackage[mathscr]{eucal}
\usepackage{times}
\def\N{\mathbb{N}}
\def\R{\mathbb{R}}

\def\C{\mathscr{C}}
\def\G{\mathscr{G}}
\def\H{\mathscr{H}}
\def\S{\mathscr{S}}

\def\M{\mathscr{M}}
\def\A{\mathscr{A}}
\def\B{\mathscr{B}}
\def\CM{\mathscr{CM}}
\def\CP{\mathscr{CP}}

\def\muhat{\widehat{\mu}}

\newtheorem{theorem}{Theorem}[section]
\newtheorem*{theorem*}{Theorem}
\def\Thm#1#2{\ifthenelse{\equal{#1}{*}}{\begin{theorem*}#2\end{theorem*}}
  {\begin{theorem}\label{T#1}#2\end{theorem}}}
\newtheorem{Atheorem}{Theorem}

\def\thm#1{Theorem~\ref{T#1}}

\newtheorem{proposition}[theorem]{Proposition}
\newtheorem*{proposition*}{Proposition}
\def\Prp#1#2{\ifthenelse{\equal{#1}{*}}{\begin{proposition*}#2\end{proposition*}}
             {\begin{proposition}\label{P#1}#2\end{proposition}}}
\def\prp#1{Proposition~\ref{P#1}}

\newtheorem{corollary}[theorem]{Corollary}
\newtheorem*{corollary*}{Corollary}
\def\Cor#1#2{\ifthenelse{\equal{#1}{*}}{\begin{corollary*}#2\end{corollary*}}
             {\begin{corollary}\label{C#1}#2\end{corollary}}}

\newtheorem{lemma}[theorem]{Lemma}
\newtheorem*{lemma*}{Lemma}
\def\Lem#1#2{\ifthenelse{\equal{#1}{*}}{\begin{lemma*}#2\end{lemma*}}
             {\begin{lemma}\label{L#1}#2\end{lemma}}}
\def\lem#1{Lemma~\ref{L#1}}

\newtheorem{remark}[theorem]{Remark}
\newtheorem*{remark*}{Remark}
\def\Rem#1#2{\ifthenelse{\equal{#1}{*}}{\begin{remark*}\rm #2\end{remark*}}
             {\begin{remark}\label{R#1}\rm #2\end{remark}}}

\def\eq#1{{\rm(\ref{E#1})}}

\def\Eq#1#2{\ifthenelse{\equal{#1}{*}}
  {\begin{equation*}\begin{aligned}#2\end{aligned}\end{equation*}}
  {\begin{equation}\begin{aligned}\label{E#1}#2\end{aligned}\end{equation}}}

\long\def\comment#1{}
\def\Det#1#2{\left|\!\begin{array}{cc}#1&#2\end{array}\!\right|}

\begin{document}
\begin{flushright}
%\textit{Submitted to: Publ.\ Math.\ Debrecen}
\end{flushright}
\vspace{5mm}

\date{\today}
\title{Characterizations of the equality of two-variable generalized quasiarithmetic means}

\author[Zs. P\'ales]{Zsolt P\'ales}
\address[Zs. P\'ales]{Institute of Mathematics, University of Debrecen, H-4002 Debrecen, Pf.\ 400, Hungary}
\email{pales@science.unideb.hu}

\author[A. Zakaria]{Amr Zakaria}
\address[A. Zakaria]{Institute of Mathematics, University of Debrecen, H-4002 Debrecen, Pf.\ 400, Hungary; Department of Mathematics, Faculty of Education, Ain Shams University, Cairo 11341, Egypt}
\email{amr.zakaria@science.unideb.hu}


%\dedicatory{}

\thanks{The research of the first author was supported by the K-134191 NKFIH Grant and the 2019-2.1.11-TÉT-2019-00049, the EFOP-3.6.1-16-2016-00022 and the EFOP-3.6.2-16-2017-00015 projects. The last two projects are co-financed by the European Union and the European Social Fund. The research of the second author was supported by the Bilateral State Scholarship of the Tempus Public Foundation of Hungary BE AK 2020-2021/157507.} 

\begin{abstract}
This paper is motivated by an astonishing result of H.\ Alzer and S.\ Ruscheweyh published in 2001 in the Proc.\ Amer.\ Math.\ Soc., which states that the intersection of the classes two-variable Gini means and Stolarsky means is equal to the class of two-variable power means. The two-variable Gini and Stolarsky means form two-parameter classes of means expressed in terms of power functions. They can naturally be generalized in terms of the so-called Bajraktarevi\'c and Cauchy means. Our aim is to show that the intersection of these two classes of functional means, under high-order differentiability assumptions, is equal to the class of two-variable quasiarithmetic means. 
\end{abstract}

\maketitle

\section{\bf Introduction}

A notion, which subsumes the concept of arithmetic, geometric and harmonic means is the concept of power means. For a real number $a$, the \textit{two-variable $a^{\mbox{\footnotesize th}}$-H\"older or $a^{\mbox{\footnotesize th}}$-power mean} $\H_a:\R_+^2\to\R$ is classically defined as
\Eq{*}{
  \H_a(x,y)
     :=\begin{cases}
       \left(\dfrac{x^a+y^a}{2}\right)^{\frac{1}{a}}
          &\mbox{ if } a\neq0,\\
       \sqrt{xy}&\mbox{ if } a=0.
       \end{cases}
}
Observe that, for $a=1$, $a=0$, and $a=-1$, the power mean $\H_a$ equals the arithmetic, geometric and harmonic means, respectively.
The theory of power means is well-developed, most of the details of their theory can be found in the monographs \cite{Bul03}, \cite{BulMitVas88}, \cite{HarLitPol34}, \cite{Mit70}, and \cite{MitPecFin93}.

The class of two-variable power means has been extended in numerous ways in the literature. One early extension was introduced by C.\ Gini \cite{Gin38} in 1938 who, for two real parameters $a,b$, defined the mean $\G_{a,b}:\R_+^2\to\R_+$ by
\Eq{Gini}{
   \G_{a,b}(x,y)
      :=\begin{cases}
         \bigg(\dfrac{x^a+y^a}{x^b+y^b}\bigg)^{\frac1{a-b}}&\mbox{if }a\neq b,\\[4mm]
         \exp\bigg(\dfrac{x^a\log(x)+y^a\log(y)}{x^a+y^a}\bigg)&\mbox{if }a=b.
        \end{cases}
}
These means are nowadays called  \textit{two-variable Gini means}. One can easily see that the two-variable H\"older means form a subclass of two-variable Gini means. Indeed, for $a\in\R$, we have $\H_a=\G_{a,0}=\G_{0,a}$ and $\H_0=\G_{a,-a}$.

Another extension of the class of two-variable power means was discovered by K.\ Stolarsky \cite{Sto75} in 1975, who, for two real parameters $a,b$, defined $\S_{a,b}:\R_+^2\to\R_+$ by
\Eq{Sto}{
   \S_{a,b}(x,y)
      :=\begin{cases}
         \bigg(\dfrac{b(x^a-y^a)}{a(x^b-y^b)}\bigg)^{\frac1{a-b}}&\mbox{if }ab(a-b)(x-y)\neq0,\\[4mm]
         \exp\bigg(-\dfrac1a+\dfrac{x^a\log(x)-y^a\log(y)}{x^a-y^a}\bigg)&\mbox{if }a=b, ab(x-y)\neq0,\\[4mm]
         \bigg(\dfrac{x^a-y^a}{a(\log(x)-\log(y))}\bigg)^{\frac1{a}}&\mbox{if }b=0, a(x-y)\neq0,\\[4mm]
         \bigg(\dfrac{x^b-y^b}{b(\log(x)-\log(y))}\bigg)^{\frac1{b}}&\mbox{if }a=0, b(x-y)\neq0,\\[4mm]
         \sqrt{xy}                                                  &\mbox{if } a=b=0,\\[4mm]
         x                                                          &\mbox{if }x=y.
        \end{cases}
}
We call these means \textit{Stolarsky means} today. We note, that these means are sometimes called \emph{extended means} (cf. \cite{LeaSho78,LeaSho83}) or difference means \cite{Pal88b,Pal92a}. One can also see that the two-variable H\"older means form a subclass of Stolarsky means. Indeed, for $a\in\R$, it is easy to see that $\H_a=\S_{2a,a}=\S_{a,2a}$, and $\H_0=\S_{a,-a}$ hold.

Therefore, the class of two-variable H\"older means is contained in the intersection of the classes of two-variable Gini and Stolarsky means. In 2001, H.\ Alzer and St. Ruscheweyh \cite{AlzRus01} established a surprising result which asserts that, instead of inclusion, the equality holds here, i.e., the class of two-variable H\"older means is equal to the intersection of the classes of two-variable Gini and Stolarsky means. 

In what follows, we will recall further important classes of two-variable functional means that extend H\"older, Gini and Stolarsky means in a natural way. These classes are the quasiarithmetic, Bajraktarevi\'c, and Cauchy means. Motivated by the above-described result of Alzer and Ruscheweyh \cite{AlzRus01}, our aim is to show that, under 8 times differentiability assumptions, the intersection of the classes of two-variable Bajraktarevi\'c and Cauchy means equals the class of two-variable quasiarithmetic means.

Throughout this paper $I$ will stand for a nonempty open real interval. In the sequel, the classes of continuous strictly monotone and continuous positive real-valued functions defined on $I$ will be denoted by $\CM(I)$ and $\CP(I)$, respectively.

Given a function $\varphi\in\CM(I)$, the \textit{two-variable quasiarithmetic mean} generated by $\varphi$ is the map $\A_\varphi:I^2\to I$ defined by
\Eq{*}{
  \A_\varphi(x,y)
      :=\varphi^{-1}\bigg(\frac{\varphi(x)+\varphi(y)}{2}\bigg).
}
The systematic treatment of these means was first given in the monograph \cite{HarLitPol34} by Hardy, Littlewood and P\'olya. The most basic problem, the characterization of the equality of these means, states that $\A_\varphi$ and $\A_\psi$ are equal to each other if and only if there exist two real constants $a\neq0$ and $b$ such that $\psi=a\varphi+b$. Observe that, by taking power functions and the logarithmic function as generators, every power mean turns out to be a quasiarithmetic mean.

In 1958, M.\ Bajraktarevi\'c \cite{Baj58}, \cite{Baj63} introduced a generalization of quasiarithmetic means essentially in the following way: For $f,g:I\to\R$ such that $f/g\in\CM(I)$ and $g\in\CP(I)$, define $\B_{f,g}:I^2\to I$ by
\Eq{*}{
  \B_{f,g}(x,y)
      :=\bigg(\frac{f}{g}\bigg)^{-1}\bigg(\frac{f(x)+f(y)}{g(x)+g(y)}\bigg).
}
The mean so defined will be called a \emph{two-variable Bajraktarevi\'c mean}. It is clear that, by taking $g=1$, the mean $\B_{f,g}$ equals $\A_f$, therefore, the class of quasiarithmetic means forms a subclass of Bajraktarevi\'c means. Assuming 6 times 
continuous differentiability, the equality problem of these means was solved by Losonczi \cite{Los99,Los06b}. A recent characterization of this equality in terms of eight equivalent conditions has been established in \cite[Theorem 15]{LosPalZak21}.

Another important generalization of quasiarithmetic means was  introduced as follows: If $f,g:I\to \R$ are continuously differentiable functions with $g'\in\CP(I)$ and $f'/g'\in\CM(I)$, then define $C_{f,g}:I^2\to I$ by
\Eq{*}{
  C_{f,g}(x,y)
      :=\begin{cases}
        \bigg(\dfrac{f'}{g'}\bigg)^{-1}
             \bigg(\dfrac{f(x)-f(y)}{g(x)-g(y)}\bigg)
            &\mbox{if }x\neq y,\\
         x  &\mbox{if }x=y.
        \end{cases}
}
The mean value property of this mean is a direct consequence of the Cauchy Mean Value Theorem, this is why this mean is called a \emph{Cauchy mean or difference mean} in the literature (cf. \ \cite{LeaSho84}, \cite{Los00a}). It is obvious that if $\varphi$ is differentiable with a nonvanishing derivative, then the mean $C_{\varphi^2,\varphi}$ equals $\A_{\varphi}$. Consequently, the class of quasiarithmetic means (with a differentiable generator) forms also a subclass of Cauchy means. Assuming 7 times continuous differentiability, the equality problem of these means was solved by Losonczi \cite{Los00a}. A characterization of this equality in terms of eight equivalent conditions has also been established in paper \cite[Theorem 16]{LosPalZak21}.

In this paper, we recall the following generalization of quasiarithmetic means, which was introduced in \cite{LosPal08} and also investigated in \cite{LosPal11a}. Given two continuous functions $f,g:I\to\R$ with $g\in\CP(I)$, $f/g\in\CM(I)$ and a probability measure $\mu$ on the Borel subsets of $[0,1]$, the two-variable mean $\M_{f,g;\mu}:I^2\to I$ is defined by
\Eq{*}{
   \M_{f,g;\mu}(x,y)
      :=\bigg(\frac{f}{g}\bigg)^{-1}\left(
            \frac{\int_{[0,1]} f\big(tx+(1-t)y\big)d\mu(t)}
                 {\int_{[0,1]} g\big(tx+(1-t)y\big)d\mu(t)}\right).
}
Means of the above form, will be called \emph{generalized quasiarithmetic means}.

In what follows, let $\delta_\tau$ denote the Dirac measure concentrated at the point $\tau\in[0,1]$. Using this notation, one can see that if $\mu=\frac{\delta_0+\delta_1}{2}$, then $\M_{f,g;\mu}$ is equal to $\B_{f,g}$ provided that $g\in\CP(I)$ and $f/g\in\CM(I)$. 

On the other hand, if $\mu$ is equal to the Lebesgue measure $\lambda$ restricted to $[0,1]$, and $g'\in\CP(I)$, $f'/g'\in\CM(I)$, then, using the Fundamental Theorem of Calculus, the mean $\M_{f',g';\mu}$ simplifies to $C_{f,g}$. Therefore, the classes of Bajraktarevi\'c and Cauchy means form subclasses of generalized quasiarithmetic means.  

The equality problem of means in various classes of two-variable means has been investigated and solved by now. We refer here to Losonczi's works \cite{Los99}, \cite{Los00a}, \cite{Los02a}, \cite{Los03a}, \cite{Los06b}, where the equality of two-variable means is characterized in various settings. A key idea in these papers, under high order differentiability assumptions, is to calculate and then to compare the partial derivatives of the means at diagonal points of the Cartesian product $I\times I$. A similar problem, the mixed equality problem of quasiarithmetic and Lagrangian means was solved by P\'ales \cite{Pal11}.

The equality problem of generalized quasiarithmetic means with the same probability measure $\mu$, i.e., the characterization of those pairs of functions $(f,g)$ and $(h,k)$ such that
\Eq{*}{
   \M_{f,g;\mu}(x,y)=\M_{h,k;\mu}(x,y)   \qquad(x,y\in I)
}
holds, was investigated and partially solved in the paper \cite{LosPalZak21}. In the particular cases $\mu=\frac{\delta_0+\delta_1}{2}$ and $\mu=\lambda$ of these results, the equality problem of two-variable Bajraktarevi\'c and the equality problem of Cauchy means was solved under $6^{\mbox{\footnotesize th}}$-order differentiability assumptions.

The aim of this paper is to study the equality problem of generalized quasiarithmetic means with the possibly different probability measures $\mu$ and $\nu$. In other words, we aim to characterize those pairs of functions $(f,g)$ and $(h,k)$ such that
\Eq{*}{
   \M_{f,g;\mu}(x,y)=\M_{h,k;\nu}(x,y)   \qquad(x,y\in I).
}
The investigation of this functional equation requires $8^{\mbox{\footnotesize th}}$-order differentiability assumptions.
Due to the complexity of this problem, we will assume that the measures are symmetric with respect to the midpoint of the interval $[0,1]$. The final main goal is to solve this equation when $\mu=\frac{\delta_0+\delta_1}{2}$ and $\nu=\lambda$. This is exactly the problem of equality of two-variable Bajraktarevi\'c means to Cauchy means. As a consequence of this result, it will follow that the intersection of these two classes of means consists of quasiarithmetic means.

\section{\bf Auxiliary results on measures}
\setcounter{equation}{0}

Given a Borel probability measure $\mu$ on the interval $[0,1]$, we define
the \emph{$k$th moment} and the \emph{$k$th centralized moment} of $\mu$ by
\Eq{*}{
  \muhat_k:=\int_0^1 t^k d\mu(t) \qquad\mbox{and}\qquad
  \mu_k:=\int_0^1 (t-\muhat_1)^k d\mu(t) \qquad(k\in\N\cup\{0\}).
}
Clearly, $\muhat_0=\mu_0=1$ and $\mu_1=0$. In view of the binomial theorem,
we easily obtain
\Eq{*}{
  \mu_k=\int_0^1 (t-\muhat_1)^k d\mu(t)
    =\int_0^1\sum_{i=0}^k(-1)^i\binom{k}{i}t^i\muhat_1^{k-i}d\mu(t)
    =\sum_{i=0}^k(-1)^k\binom{k}{i}\muhat_i\muhat_1^{k-i}
   \qquad (k\in\N)
}
and
\Eq{muhat}{
  \muhat_k=\int_0^1 \big((t-\muhat_1)+\muhat_1\big)^k d\mu(t)
    =\int_0^1\sum_{i=0}^k\binom{k}{i}(t-\muhat_1)^i\muhat_1^{k-i}d\mu(t)
    =\sum_{i=0}^k\binom{k}{i}\mu_i\muhat_1^{k-i}
   \qquad (k\in\N).
}
It is obvious that $\mu_{2k}\geq0$ for all $k\in\N$. If , for some $k\in\N$ the equality $\mu_{2k}=0$ holds then $\mu$ has to be the Dirac measure $\delta_{\muhat_1}$. (In the sequel, $\delta_\tau$ will denote the Dirac measure concentrated at the point $\tau\in[0,1]$.) On the other hand, the odd-order centralized moments can be zero. One can prove that $\mu_{2k-1}=0$ holds for all $k\in\N$ if and only if $\mu$ is symmetric with respect to its first moment $\muhat_1$, i.e., if $\mu(A)=\mu\big((2\mu_1-A)\cap[0,1]\big)$ for all Borel sets $A\subseteq[0,1]$.

In what follows we define a two-parameter class of Borel measures which will be instrumental for our investigations.

For two given positive parameters $\ell$ and $p$, we define a probability measure $\pi:=\pi(\ell,p)$ via the following equalities for its moments
\Eq{*}{
  \pi_0:=1, \quad\hat{\pi}_1:=\tfrac12, \quad 
  \pi_{2n-1}:=0, \quad 
  \pi_{2n}:=\frac{(2n)!}{n!}\ell^n p^{\langle n\rangle}
  \qquad(n\in\N),
}
where, the modified power $p^{\langle n\rangle}$ is defined by  
\Eq{*}{
p^{\langle n\rangle}:=\prod\limits_{i=0}^{n-1}\frac{p}{1+ip} \qquad(n\in\N).
}
Clearly, $p^{\langle 1\rangle}=p^1=p$ and the recursive formula $p^{\langle n+1\rangle}=\frac{p}{1+np}\cdot p^{\langle n\rangle}$ holds.

It follows from equality \eq{muhat} that the moments of the measure $\pi(\ell,p)$ are completely determined and, therefore, by a classical result related to the Hausdorff moment problem, the measure $\pi(\ell,p)$ is uniquely determined, however, it may not exist for every $\ell,p>0$. It also follows that $\pi(\ell,p)$ has to be a symmetric measure with respect to the point $\tfrac12$. The set of those parameters $(\ell,p)$ for which the probability measure $\pi(\ell,p)$ exists will be denoted by $\Pi$.

\Lem{pi}{For the parameter set $\Pi$, we have the following inclusion
\Eq{*}{
  \Pi\subseteq \,]0,\tfrac1{16}]\times\,]0,2].
}}

\begin{proof} Let $(\ell,p)\in\Pi$. Then, due to the estimate $|t-\tfrac12|\leq\tfrac12$, we get
\Eq{*}{
\frac{(2n)!}{n!}\ell^n p^{\langle n\rangle}=\pi_{2n}(\ell,p)
=\int_{[0,1]}(t-\tfrac12)^{2n}d\pi(\ell,p)(t)
\leq \int_{[0,1]}\frac{1}{2^{2n}}d\pi(\ell,p)(t)=\frac{1}{4^n}
}
for all $n\in\N$. This inequality yields
\Eq{ub}{
4\ell
\leq \sqrt[n]{\frac{n!}{(2n)!p^{\langle n\rangle}}}.
}
In order to compute the limit of the right hand side, we shall use the multiplicative version of the Cesaro--Stolz theorem. Denote by $a_n$ the expression under the $n$th root on the right hand side of \eq{ub}. Then, according to this classical result, we get 
\Eq{*}{
 4\ell
 \leq \lim_{n\to\infty} \sqrt[n]{a_n}
 =\lim_{n\to\infty}\frac{a_{n+1}}{a_n}
 =\lim_{n\to\infty}\frac{(n+1)(1+np)}{(2n+2)(2n+1)p} =\frac{1}{4},
}
which implies the inequality $\ell\leq\tfrac1{16}$.

To prove the inequality $p\leq 2$, we apply the Cauchy--Bunyakovski--Schwartz inequality:
\Eq{*}{
  \big(\pi_{2}(\ell,p)\big)^2
  &=\bigg(\int_{[0,1]}(t-\tfrac12)^{2}\cdot 1\,d\pi(\ell,p)(t)\bigg)^2\\
  &\leq \bigg(\int_{[0,1]}(t-\tfrac12)^{4}d\pi(\ell,p)(t)\bigg)
  \bigg(\int_{[0,1]}1 \,d\pi(\ell,p)(t)\bigg)
  =\pi_{4}(\ell,p),
}
which reduces to
\Eq{*}{
  \big(2\ell p)^2\leq\frac{12(\ell p)^2}{1+p}.
}
This simplifies to the inequality $p\leq 2$.
\end{proof}

\Lem{mn}{Let $\tau\in[0,\tfrac12[\,$. Then
\Eq{*}{
  \tfrac12(\delta_\tau+\delta_{1-\tau})
  =\pi\big((\tfrac{1-2\tau}4)^2,2\big) 
  \qquad\mbox{and}\qquad 
  \tfrac1{1-2\tau}\lambda|_{[\tau,1-\tau]}
  =\pi\big((\tfrac{1-2\tau}4)^2,\tfrac23\big),
}
where $\lambda$ denotes the standard Lebesgue measure restricted to $[0,1]$.}

\begin{proof}
Denote $\mu:=\tfrac12(\delta_\tau+\delta_{1-\tau})$. Then, $\mu_0=1$, $\hat{\mu}_1=\tfrac12$ and, for $n\in\N$, we have  
\Eq{*}{
  \mu_{n}=\int_{[0,1]}(t-\tfrac12)^{n}d\mu(t)
  =\frac{(\tau-\tfrac12)^{n}+(1-\tau-\tfrac12)^{n}}{2}
  =\begin{cases}
   (\tfrac12-\tau)^{n} &\mbox{ if $n$ is even},\\ 
   0 &\mbox{ if $n$ is odd}.
   \end{cases}
}
On the other hand, if $p=2$, $\ell=(\tfrac{1-2\tau}4)^2$, then, for all $n\in\N$, $\pi_{2n-1}((\tfrac{1-2\tau}4)^2,2)=0=\mu_{2n-1}$ and
\Eq{*}{
  \pi_{2n}((\tfrac{1-2\tau}4)^2,2)
  =\frac{(2n)!}{n!}\ell^n p^{\langle n\rangle}
  =\frac{(2n)!\cdot (2\ell)^n}{n!\cdot 1\cdot3\cdots(2n-1)}
  =(4\ell)^n=(\tfrac12-\tau)^{2n}=\mu_{2n}.
}
Therefore, all the moments of $\mu$ and $\pi((\tfrac{1-2\tau}4)^2,2)$ are the same, which  proves that $\mu=\pi((\tfrac{1-2\tau}4)^2,2)$.

For the second assertion, denote $\nu:=\tfrac1{1-2\tau}\lambda|_{[\tau,1-\tau]}$. Then, $\mu_0=1$, $\hat{\mu}_1=\tfrac12$ and, for $n\in\N$, we have  
\Eq{*}{
  \nu_{n}=\int_{[0,1]}(t-\tfrac12)^{n}d\nu(t)
  &=\frac1{1-2\tau}\int_{[\tau,1-\tau]}(t-\tfrac12)^{n}dt\
  =\frac1{1-2\tau}\bigg[\frac{(t-\tfrac12)^{n+1}}{n+1}\bigg]_{t=\tau}^{t=1-\tau}\\
  &=\frac{(1-\tau-\tfrac12)^{n+1}-(\tau-\tfrac12)^{n+1}}{(1-2\tau)(n+1)}
  =\begin{cases}
   \frac{(\tfrac12-\tau)^{n}}{n+1} &\mbox{ if $n$ is even},\\ 
   0 &\mbox{ if $n$ is odd}.
   \end{cases}
}
On the other hand, if $p=\tfrac23$, $\ell=(\tfrac{1-2\tau}4)^2$, then, for all $n\in\N$,
\Eq{*}{
  \pi_{2n}((\tfrac{1-2\tau}4)^2,\tfrac23)
  =\frac{(2n)!}{n!}\ell^n p^{\langle n\rangle}
  =\frac{(2n)!\cdot (2\ell)^n}{n!\cdot 3\cdot5\cdots(2n+1)}
  =\frac{(4\ell)^n}{2n+1}=\frac{(\tfrac12-\tau)^{2n}}{2n+1}=\nu_{2n}.
}
Therefore, all the moments of $\nu$ and $\pi((\tfrac{1-2\tau}4)^2,\tfrac23)$ are the same, which  proves that $\nu=\pi((\tfrac{1-2\tau}4)^2,\tfrac23)$.
\end{proof}

\comment{Consider now the family $\pi(\ell,1)$. Then
\Eq{*}{
  \pi_{2n}(\ell,1)=\frac{(2n)!\ell^n}{n!^2}=\binom{2n}{n}\ell^n.
}
Then, if $P(t)=\sum_{k=0}^n\frac{P^{(k)}(\tfrac12)}{k!}{�}(t-\tfrac12)^k$, then
\Eq{*}{
  \int_{[0,1]}P(t)d\pi(\ell,1)(t)
  =\sum_{k=0}^{n/2}\binom{2k}{k}\frac{P^{(2k)}(\tfrac12)}{(2k)!}\ell^{2k}
  =\sum_{k=0}^{n/2}\frac{P^{(2k)}(\tfrac12)}{(k!)^2}\ell^{2k}.
}}

\section{\bf Auxiliary results on generalized quasiarithmetic means}
\setcounter{equation}{0}

In order to describe the regularity conditions related to the two functions $f,g$ generating the mean $\M_{f,g;\mu}$, we introduce some notations. The class $\C_0(I)$ consists of all those pairs $(f,g)$ of continuous functions $f,g:I\to\R$ such that $g\in\CP(I)$ and $f/g\in\CM(I)$. For $n\in\N$, we say that the pair $(f,g)$ is in the class $\C_n(I)$ if $f,g$ are $n$-times continuously differentiable functions such that $g\in\CP(I)$ and the function $f'g-fg'$ does not vanish anywhere on $I$. This latter
condition implies that $f/g$ is strictly monotone, i.e., $f/g\in\CM(I)$. Therefore, $\C_0(I)\supseteq\C_1(I)\supseteq \C_2(I)\supseteq\cdots$.

For $(f,g)\in\C_2(I)$, we also introduce the notation
\Eq{*}{
  \Phi_{f,g}:=\frac{W_{f,g}^{2,0}}{W_{f,g}^{1,0}}\qquad\mbox{and}\qquad
  \Psi_{f,g}:=-\frac{W_{f,g}^{2,1}}{W_{f,g}^{1,0}},
}
where the $(i,j)$-order Wronskian operator $W_{f,g}^{i,j}$ is defined in terms of the $i$th and $j$th derivatives by
\Eq{*}{
  W_{f,g}^{i,j}:=\begin{vmatrix}f^{(i)}&f^{(j)}\\g^{(i)}&g^{(j)}\end{vmatrix}.
}
We recall now a result from the paper \cite[Lemma 1]{LosPalZak21} which establishes a formula for the $(i,j)$-order Wronskian in terms of the functions $\Phi_{f,g}$ and $\Psi_{f,g}$. 

\Lem{0}{Let $(f,g)\in\C_n(I)$, where $n\geq2$ and define two sequences
$(\varphi_i)$ and $(\psi_i)$ by the recursions
\Eq{Rec}{
      \varphi_0&:=0,&\qquad
  \varphi_{i+1}&:=\varphi'_i+\varphi_i\Phi_{f,g}+\psi_i&
      \qquad (i&\in\{0,\dots,n-1\}),\\
         \psi_0&:=1,&
     \psi_{i+1}&:=\varphi_i\Psi_{f,g}+\psi'_i &
      \qquad (i&\in\{0,\dots,n-1\}).
}
Then
\Eq{W}{
  W_{f,g}^{i,j}=\begin{vmatrix}\varphi_i&\varphi_j\\\psi_i&\psi_j\end{vmatrix}
               \cdot W_{f,g}^{1,0}\qquad (i,j\in\{0,\dots,n\}).
}
In particular, 
\Eq{W01}{
  W_{f,g}^{i,0}=\varphi_i W_{f,g}^{1,0},\qquad
  W_{f,g}^{i,1}=-\psi_i W_{f,g}^{1,0}, \qquad
  W_{f,g}^{i,2}=(\varphi_i\Psi_{f,g}-\psi_i\Phi_{f,g}) W_{f,g}^{1,0}
  \qquad (i\in\{0,\dots,n\}).
}}

For small $i$, the first few members of the sequences $(\varphi_i)$ and $(\psi_i)$ are as follows:
\Eq{*}{
   \varphi_1&=1, &\qquad \varphi_2&=\Phi_{f,g}, 
      &\qquad \varphi_3&=\Phi'_{f,g}+\Phi^2_{f,g}+\Psi_{f,g},\\
   \psi_1&=0, &\psi_2&=\Psi_{f,g},
      & \psi_3&=\Phi_{f,g}\Psi_{f,g}+\Psi'_{f,g}.
}
The subsequent elements can easily be computed by the recursion \eq{Rec}. In the sequel we shall need the following consequence of the recursion \eq{Rec}.

\Lem{DRec}{Under the same assumptions as in \lem{0}, for the sequences $(\varphi_i)$ and $(\psi_i)$, we have
\Eq{DRec}{
  \varphi_{i+2}
  &=\varphi''_i+2\varphi'_i\Phi_{f,g}+\varphi_i\varphi_3
  +2\psi'_i+\psi_i\Phi_{f,g}&\qquad (i&\in\{0,\dots,n-2\}),\\
  \psi_{i+2}
  &=2\varphi'_i\Psi_{f,g}+\varphi_i\psi_3+\psi''_{i}+\psi_i\Psi_{f,g}&
  \qquad (i&\in\{0,\dots,n-2\}).
}}

\begin{proof} Applying the recursion \eq{Rec} twice, for $i\in\{0,\dots,n-2\}$, we get
 \Eq{*}{
  \varphi_{i+2}
  &=\varphi'_{i+1}+\varphi_{i+1}\Phi_{f,g}+\psi_{i+1}\\
  &=(\varphi'_i+\varphi_i\Phi_{f,g}+\psi_i)'+(\varphi'_i+\varphi_i\Phi_{f,g}+\psi_i)\Phi_{f,g}+\varphi_i\Psi_{f,g}+\psi'_i,\\
  &=\varphi''_i+2\varphi'_i\Phi_{f,g}+\varphi_i\varphi_3
  +2\psi'_i+\psi_i\Phi_{f,g},\\
  \psi_{i+2}
  &=\varphi_{i+1}\Psi_{f,g}+\psi'_{i+1}\\
  &=(\varphi'_i+\varphi_i\Phi_{f,g}+\psi_i)\Psi_{f,g}
  +(\varphi_{i}\Psi_{f,g}+\psi'_{i})'\\
  &=2\varphi'_i\Psi_{f,g}+\varphi_i\psi_3
  +\psi''_{i}+\psi_i\Psi_{f,g},
}
which completes the proof of \eq{DRec}.
\end{proof}

The remaining auxiliary results of this section were obtained in \cite{LosPal08} and \cite{PalZak20a}. The next assertion characterizes the mean $\M_{f,g;\mu}$ via an implicit equation.

\Lem{1}{{\rm(\cite[Lemma 1]{LosPal08}, \cite[Lemma 1.1]{PalZak20a})} Let $(f,g)\in\C_0(I)$ and $\mu$ be a Borel probability measure on $[0,1]$. Then for all $x,y\in I$, the value $z=\M_{f,g;\mu}(x,y)$ is the unique solution of the equation
\Eq{1}{
  \int_{[0,1]} \begin{vmatrix}f\big(tx+(1-t)y\big) & f(z)\\
            g\big(tx+(1-t)y\big) & g(z)\end{vmatrix}d\mu(t)=0.
}}

We say that \textit{two pairs of functions $(f,g),(h,k)\in\C_0(I)$ are equivalent}, denoted by $(f,g)\sim(h,k)$, if there exist four real constants $a,b,c,d$ with $ad\neq bc$ such that 
\Eq{*}{
h=af+bg\qquad\mbox{and}\qquad k=cf+dg.
}

The property of equivalence in the class $\C_2(I)$ is completely characterized by the following result.

\Lem{3}{{\rm(\cite[Theorem 2.1]{PalZak20a})} Let $(f,g),(h,k)\in\C_2(I)$. Then $(f,g)\sim(h,k)$ holds if and only if
\Eq{*}{
  \Phi_{f,g}=\Phi_{h,k}\qquad\mbox{and}\qquad \Psi_{f,g}=\Psi_{h,k}.
}}

As a consequence \lem{1}, the next lemma shows that the equivalent pairs of generating functions determine identical means.

\Lem{2}{{\rm(\cite{LosPal08}, \cite{PalZak20a})} Let $(f,g),(h,k)\in\C_0(I)$ and $\mu$ be a Borel probability measure on $[0,1]$. Assume that $(f,g)\sim(h,k)$. Then $\M_{f,g;\mu}=\M_{h,k;\mu}$.}

For a real parameter $a\in\R$, introduce the sine and cosine type functions $S_a,C_a:\R\to\R$ by
\Eq{*}{
  S_a(x):=\begin{cases}
           \sin(\sqrt{-a}x) & \mbox{ if } a<0, \\
           x & \mbox{ if } a=0, \\
           \sinh(\sqrt{a}x) & \mbox{ if } a>0, 
         \end{cases}\qquad\mbox{and}\qquad
  C_a(x):=\begin{cases}
           \cos(\sqrt{-a}x) & \mbox{ if } a<0, \\
           1 & \mbox{ if } a=0, \\
           \cosh(\sqrt{a}x) & \mbox{ if } a>0. 
        \end{cases}
}
It is easily seen that the functions $S_a$ and $C_a$ form the fundamental system of solutions for the second-order homogeneous linear differential equation $h''=ah$.

\section{\bf Higher-order directional derivatives of generalized quasiarithmetic means}
\setcounter{equation}{0}

We first recall the result stated in \cite[Lemma5]{LosPalZak21}.

\Lem{n}{Let $n\in\N$, $(f,g)\in\C_n(I)$ and $\mu$ be a Borel probability measure on $[0,1]$. Then $\M_{f,g;\mu}$ is $n$-times continuously differentiable on $I\times I$.}

In what follows, we deduce explicit formulae for the high-order directional derivatives of $\M_{f,g;\mu}$ at the diagonal points of the Cartesian product $I\times I$ provided that $\mu$ is a symmetric probability measure. Due to the symmetry of $\mu$, we can see that $\hat{\mu}_1=\frac12$ and $\mu_{2n-1}=0$ for all $n\in\N$ and that $\M_{f,g;\mu}$ is a symmetric two-variable mean. 

Given $(f,g)\in\C_0(I)$ and a fixed element $x\in I$, define the function $m_x=m_{x;f,g;\mu}$ in the symmetric neighborhood $U_x:=2(I-x)\cap2(x-I)$ of zero by
\Eq{m}{
  m_x(u):=m_{x;f,g;\mu}(u):=\M_{f,g;\mu}\big(x+\tfrac12 u,x-\tfrac12 u\big).
}
The symmetry of $\mu$ yields that $m_x:U_x\to\R$ is even, that is, $m_x(u)=m_x(-u)$ holds for all $u\in U_x$. 

In the subsequent results, the derivatives of any real function $h$ up to the order 8 will be denoted by $h^{(i)}$ for all $i\in\{0,\dots,8\}$. Alternatively, $h^{(0)},h^{(1)},h^{(2)}$ and $h^{(3)}$, will be denoted by $h,h',h''$ and $h'''$, respectively.

\Prp{der}{Let $n\in\N$, $(f,g)\in\C_n(I)$, and $\mu$ be a symmetric Borel probability measure on $[0,1]$. Then, for fixed $x\in I$, the function $m_x$ defined by \eq{m} is $n$-times differentiable at the origin. Furthermore, $m_x(0)=x$ and, in the cases $n\in\{1,3,5,7\}$, the equality $m_x^{(n)}(0)=0$ holds and, in the cases $n\in\{2,4,6,8\}$, we have
\Eq{*}{
    m_x''(0)&=\mu_2\varphi_2(x),\\
    m_x^{(4)}(0)&=\mu_4\varphi_4(x)-3\mu_2^2(\varphi_2^3+2\psi_2\varphi_2)(x),\\
  m_x^{(6)}(0)&=\mu_6\varphi_6(x)-15\mu_4\mu_2(\varphi_4(\varphi_2^2+\psi_2)+\varphi_2\psi_4)(x)-15\mu_2^3\varphi_2(\varphi_3\varphi_2^2
  -3(\varphi_2^2+\psi_2)(\varphi_2^2+2\psi_2))(x),\\
  m_x^{(8)}(0)&=\mu_8\varphi_8(x)-28\mu_6\mu_2(\varphi_6(\varphi_2^2+\psi_2)
  +\varphi_2\psi_6)(x)-35\mu_4^2(\varphi_4^2\varphi_2+2\varphi_4\psi_4)(x)\\
  &\quad+210\mu_4\mu_2^2(\varphi_4(3\varphi_2^4+\varphi_2^2(7\psi_2-\varphi_3)+2\psi_2^2)+2\varphi_2\psi_4(\varphi_2^2+2\psi_2))(x)\\
  &\quad-105\mu_2^4(\varphi_4\varphi_2^4+15\varphi_2^7+2\varphi_2^3(5\varphi_2^2+6\psi_2)(6\psi_2-\varphi_3)+4\varphi_2(6\psi_2^3-\varphi_2^3\psi_3))(x),
}
where the functions $\varphi_i=\varphi_{i;f,g}$ ($i\in\{2,3,4,6,8\}$) and $\psi_i=\psi_{i;f,g}$ ($i\in\{2,3,4,6\}$) are defined by the recursion \eq{Rec}.}

\begin{proof} Let $x\in I$ be fixed. The equality $m_x(0)=x$ follows from the definition of $m_x$. The $n$-times differentiability of $m_x$ at $u=0$ is a consequence of \lem{n}. By the symmetry of the mean $\M_{f,g;\mu}$, $m_x:U_x\to\R$ is even. This implies that all odd-order derivatives of $m_x$ vanish at the point $u=0$.

For brevity, let $F$ denote the vector valued function $\Big(\!\!\begin{array}{c}f\\g\end{array}\!\!\Big):I\to\R^2$. Then, for $u\in U_x$, \eq{1} can be rewritten as
\Eq{*}{
  \int_{[0,1]}\Det{F(x+(t-\frac{1}{2})u)}{(F\circ m_x)(u)}d\mu(t)=0.
}
Differentiating this equality $n$-times with respect to the variable $u$ and using Leibniz's rule, we obtain
\Eq{*}{
  \sum_{i=0}^n\binom{n}{i}\int_{[0,1]}
  \Det{F^{(i)}(x+(t-\frac{1}{2})u)}{(F\circ m_x)^{(n-i)}(u)}
  (t-\tfrac{1}{2})^id\mu(t)=0.
}
Now the substitution $u=0$ gives
\Eq{*}{
  \sum_{i=0}^n\binom{n}{i}
  \mu_i\Det{F^{(i)}(x)}{(F\circ m_x)^{(n-i)}(0)}=0.
}
Using that the odd-order centralized moments of $\mu$ are equal to zero, this equality finally simplifies to
\Eq{der}{
  \sum_{i=0}^{\big\lfloor \tfrac{n}{2}\big\rfloor}
  \binom{n}{2i}\mu_{2i}
  \Det{F^{(2i)}(x)}{(F\circ m_x)^{(n-2i)}(0)}=0.
}
The evenness of $m_x$ implies that $F\circ m_x$ is also even on $U_x$ and hence all of its (existing) odd-order derivatives vanish at $u=0$. This shows that \eq{der} is nontrivial only when $n$ is even. To elaborate the condition \eq{der} in the cases $n\in\{2,4,6,8\}$, we shall use Fa\'a di Bruno's formula to compute $(F\circ m_x)^{(N)}$ as follows
\Eq{Faa}{
(F\circ m_x)^{(N)}(0):=
\sum_{k=1}^N (F^{(k)}\circ m_x)\cdot \B_{N,k}\Big(m_x'(0),m_x''(0),\dots,m_x^{(N-k+1)}(0)\Big).
}
Here $\B_{N,k}:\R^{N-k+1}\to\R$ is the incomplete Bell polynomial, which is defined by the recursive formula
\Eq{*}{
\B_{N,k}(x_1,\dots,x_{N-k+1}):=\sum_{j=1}^{N-k+1}\binom{N-1}{j-1}x_j\B_{N-j,k-1}(x_1,\dots,x_{N-j-k+2}),
}
such that $\B_{0,0}=1$ and $\B_{N,0}=0=\B_{0,k}$ for all $k,N\in\N$. When applying the formula \eq{Faa}, we have that the odd-order derivatives of $m_x$ are zero at $u=0$, therefore, we need to compute the Bell polynomials only for such arguments where the odd coordinates are equal to zero. Now, an easy calculation shows that, for all $N\in\{1,\dots,8\}$, $k\in\{1,\dots,N\}$ and argument $\pmb{x}_i:=(x_1,x_2,\dots,x_i)$ where $x_{2j-1}=0$ for all $j\in\{1,\dots,\lfloor \tfrac{i+1}{2}\rfloor\}$,
\Eq{*}{
\B_{N,k}(\pmb{x}_{N-k+1})=0
\qquad\mbox{if either }N \mbox{ is odd \& } k\in\{1,\dots,N\}
\quad\mbox{or }N\mbox{ is even \& } k\in\{\tfrac{N+2}{2},\dots,N\},
}
and, for the remaining cases, we have
\Eq{*}{
\B_{2,1}(\pmb{x}_2)&=x_2,\\
\B_{4,1}(\pmb{x}_4)&=x_4,&
\B_{4,2}(\pmb{x}_3)&=3x_2^2,& \\
\B_{6,1}(\pmb{x}_6)&=x_6,&
\B_{6,2}(\pmb{x}_5)&=15x_2x_4,& 
\B_{6,3}(\pmb{x}_4)&=15x_2^3,\\
\B_{8,1}(\pmb{x}_8)&=x_8,&\qquad
\B_{8,2}(\pmb{x}_7)&=28x_2x_6+35x_4^2,&\qquad
\B_{8,3}(\pmb{x}_6)&=210x_2^2x_4,&\qquad
\B_{8,4}(\pmb{x}_5)&=105x_2^4.
}
These identities, together with formula \eq{Faa}, imply that
\Eq{ids}{
  (F\circ m_x)(0)&=F(x),\\
  (F\circ m_x)''(0)&=F'(x)m_x''(0),\\
  (F\circ m_x)^{(4)}(0)&=F'(x)m_x^{(4)}(0)+3F''(x)(m_x''(0))^2,\\
  (F\circ m_x)^{(6)}(0)&=F'(x)m_x^{(6)}(0)
          +15F''(x)m_x''(0)m_x^{(4)}(0)+15F'''(x)(m_x''(0))^3,\\
  (F\circ m_x)^{(8)}(0)&=F'(x)m_x^{(8)}(0)
  +F''(x)(28m_x''(0)m_x^{(6)}(0)+35(m_x^{(4)}(0))^2)\\
  &\quad+210F'''(x)(m_x''(0))^2m_x^{(4)}(0)
  +105F^{(4)}(x)(m_x''(0))^4.
}

From now on, we shall use the notations $\varphi_i$ and $\psi_i$ which were introduced in \lem{0}.
In the case $n=2$, using $\mu_0=1$, equation \eq{der} yields
\Eq{*}{
\Det{F(x)}{(F\circ m_x)''(0)}
+\mu_2\Det{F''(x)}{(F\circ m_x)(0)}=0,
}
which, in view of the first identity in \eq{ids}, reduces to the equality
\Eq{*}{
  -W_{f,g}^{1,0}(x)m_x''(0)+\mu_2W_{f,g}^{2,0}(x)=0.
}
This proves $m_x''(0)=\mu_2\varphi_2(x)$.

In the case $n=4$, using $\mu_0=1$, equation \eq{der} gives
\Eq{*}{
\Det{F(x)}{(F\circ m_x)^{(4)}(0)}
+6\mu_2\Det{F''(x)}{(F\circ m_x)''(0)}
+\mu_4\Det{F^{(4)}(x)}{(F\circ m_x)(0)}=0,
}
which, by the first and second identities in \eq{ids}, implies
\Eq{*}{
  -W_{f,g}^{1,0}(x)m_x^{(4)}(0)-3W_{f,g}^{2,0}(x)(m_x''(0))^2
    +6\mu_2W_{f,g}^{2,1}(x)m_x''(0)+\mu_4W_{f,g}^{4,0}(x)=0.
}
Thus, applying \eq{W01} and the formula for $m_x''(0)$, we obtain
\Eq{*}{
  m_x^{(4)}(0)
    &=-3\varphi_2(x)(m_x''(0))^2-6\mu_2\psi_2(x)m_x''(0)+\mu_4\varphi_4(x)\\
    &=-3\mu_2^2\varphi_2^3(x)-6\mu_2^2\psi_2(x)\varphi_2(x)+\mu_4\varphi_4(x).
}

In the case $n=6$, using $\mu_0=1$, equation \eq{der} results
\Eq{*}{
\Det{F(x)}{(F\circ m_x)^{(6)}(0)}
  &+15\mu_2\Det{F''(x)}{(F\circ m_x)^{(4)}(0)}\\
  &+15\mu_4\Det{F^{(4)}(x)}{(F\circ m_x)''(0)}
  +\mu_6\Det{F^{(6)}(x)}{(F\circ m_x)(0)}=0,
}
which, using the identities \eq{ids}, implies
\Eq{*}{
  -W_{f,g}^{1,0}(x)m_x^{(6)}(0)
  &-15W_{f,g}^{2,0}(x)m_x''(0)m_x^{(4)}(0)
  -15W_{f,g}^{3,0}(x)(m_x''(0))^3\\
  &+15\mu_2W_{f,g}^{2,1}(x)m_x^{(4)}(0)
  +15\mu_4W_{f,g}^{4,1}(x)m_x''(0)+\mu_6W_{f,g}^{6,0}(x)=0.
}
Thus, applying \eq{W01}, we get
\Eq{*}{
  m_x^{(6)}(0)
   =&-15\varphi_2(x)m_x''(0)m_x^{(4)}(0)-15\varphi_3(x)(m_x''(0))^3 
   \\&-15\mu_2\psi_2(x)m_x^{(4)}(0)
      -15\mu_4\psi_4(x)m_x''(0)+\mu_6\varphi_6(x).
}
Using the explicit formulae for $m_x''(0)$ and $m_x^{(4)}(0)$, this formula simplifies to the statement.

Finally, in the case $n=8$, equation \eq{der} reduces to the equality
\Eq{*}{
\Det{F(x)}{(F\circ m_x)^{(8)}(0)}
  &+28\mu_2\Det{F''(x)}{(F\circ m_x)^{(6)}(0)}  
  +70\mu_4\Det{F^{(4)}(x)}{(F\circ m_x)^{(4)}(0)}\\
  &+28\mu_6\Det{F^{(6)}(x)}{(F\circ m_x)''(0)}
  +\mu_8\Det{F^{(8)}(x)}{(F\circ m_x)(0)}=0.
}
Now, applying the identities \eq{ids}, gives
\Eq{*}{
-&W_{f,g}^{1,0}(x)m_x^{(8)}(0)
  -W_{f,g}^{2,0}(x)(28m_x''(0)m_x^{(6)}(0)+35(m_x^{(4)}(0))^2)
  -210W_{f,g}^{3,0}(x)(m_x''(0))^2m_x^{(4)}(0)\\
  &-105W_{f,g}^{4,0}(x)(m_x''(0))^4
  +28\mu_2W_{f,g}^{2,1}(x)m_x^{(6)}(0)
  -420\mu_2W_{f,g}^{3,2}(x)(m_x''(0))^3\\
  &+70\mu_4W_{f,g}^{4,1}(x)m_x^{(4)}(0)
  +210\mu_4W_{f,g}^{4,2}(x)(m_x''(0))^2
  +28\mu_6W_{f,g}^{6,1}(x)m_x''(0)
  +\mu_8W_{f,g}^{8,0}(x)=0.
}
Then, using \eq{W} and \eq{ids}, we get
\Eq{*}{
m_x^{(8)}(0)=&-\varphi_2(x)(28m_x''(0)m_x^{(6)}(0)+35(m_x^{(4)}(0))^2)
  -210\varphi_3(x)(m_x''(0))^2m_x^{(4)}(0)\\
  &-105\varphi_4(x)(m_x''(0))^4 
  -28\mu_2\psi_2(x)m_x^{(6)}(0)
  -420\mu_2(\varphi_3\psi_2-\psi_3\varphi_2)(x)(m_x''(0))^3\\
  &-70\mu_4\psi_4(x)m_x^{(4)}(0)
  +210\mu_4(\varphi_4\psi_2-\psi_4\varphi_2)(x)(m_x''(0))^2
  -28\mu_6\psi_6(x)m_x''(0)
  +\mu_8\varphi_8(x).
}
This, together with the formulae for $m_x''(0),m_x^{(4)}(0)$, and $m_x^{(6)}(0)$, implies the required result.
\end{proof}

\section{\bf Necessary and sufficient conditions for the equality of generalized quasiarithmetic means}
\setcounter{equation}{0}

In what follows, given $(f,g),(h,k)\in\C_0(I)$ and a probability measure $\mu$ on $[0,1]$, we say that $\M_{f,g;\mu}$ equals $\M_{h,k;\mu}$ if they coincide at every point of $I^2$. We say that these two means are \emph{equal near the diagonal $\Delta(I):=\{(x,x)\mid x\in I\}$ of $I^2$}, if there exists an open set $U\subseteq I^2$ containing a dense subset $D$ of $\Delta(I)$ such that the two means are equal at every point of $U$.

\Lem{N0}{{\rm(\cite[Lemma 7]{LosPalZak21})}
Let $\mu$ and $\nu$ be symmetric Borel probability measures on $[0,1]$, let $n\in\N$ and let $(f,g),(h,k)\in\C_{2n}(I)$. If $\M_{f,g;\mu}$ equals $\M_{h,k;\nu}$ near the diagonal of $I^2$, then, for all $i\in\{1,\dots,n\}$ and $x\in I$,
\Eq{mk}{
m_{x;f,g;\mu}^{(2i)}(0)=m_{x;h,k;\nu}^{(2i)}(0).
}
}

In the following result, we consider first the particular case when $\mu=\nu$ and $\Psi_{f,g}=\Psi_{h,k}$, and then we characterize the equality of the means $\M_{f,g;\mu}$ and $\M_{h,k;\nu}$.

\Thm{N1.5}{Let $\mu$ be a symmetric Borel probability measure on $[0,1]$ with $\mu_2\neq0$. Let $(f,g),(h,k)\in\C_2(I)$ such that $\Psi_{f,g}=\Psi_{h,k}$. Then the following assertions are equivalent:
\begin{enumerate}[(i)]
 \item The means $\M_{f,g;\mu}$ and $\M_{h,k;\mu}$ are equal on $I^2$.
 \item The means $\M_{f,g;\mu}$ and $\M_{h,k;\mu}$ are equal near the diagonal of  $I^2$.
 \item For all $x\in I$, the identity $m_{x;f,g;\mu}''(0)=m_{x;h,k;\mu}''(0)$ holds.
 \item The equality $\Phi_{f,g}=\Phi_{h,k}$ holds on $I$.
 \item $(f,g)\sim(h,k)$ holds.
\end{enumerate}
}

\begin{proof}
The implication (i)$\Rightarrow$(ii) is clear. The implication (ii)$\Rightarrow$(iii) is a consequence of \lem{N0}. If (iii) holds, then, according to \prp{der}, we get that $\mu_2\Phi_{f,g}=\mu_2\Phi_{h,k}$, which reduces to the equality $\Phi_{f,g}=\Phi_{h,k}$ proving (iv). Finally, \lem{3} and \lem{2} imply the assertions (iv)$\Rightarrow$(v) and (v)$\Rightarrow$(i), respectively.
\end{proof}

Due to the symmetry of the measure $\mu$, the third-order necessary condition $m_{x;f,g;\mu}'''(0)=m_{x;h,k;\mu}'''(0)$ does not imply the equality $\Psi_{f,g}=\Psi_{h,k}$ (as it is obtained in \cite{LosPalZak21}). Therefore, in the sequel, we want to consider problems where either the measures $\mu$ and $\nu$ are different or the equality $\Psi_{f,g}=\Psi_{h,k}$ cannot be established. In view of \prp{der} and \lem{0}, without any further assumption on the measures $\mu$ and $\nu$, we can easily see that the $i$th equation in \eq{mk} is a nonlinear differential equation of order $2(i-1)$ for the unknown functions $\Phi_{f,g}$, $\Psi_{f,g}$, $\Phi_{h,k}$, and $\Psi_{h,k}$. Therefore, for the solution of the equality problem of the means, it seems to be enough to take these equations for $i\in\{1,2,3,4\}$ and solve the system of differential equations so obtained. However, the integration of this system in its full generality seems to be hopeless. Therefore, we take $\mu$ and $\nu$ from the two-parameter family  of measures $\{\pi(\ell,p):(\ell,p)\in\Pi\}$ introduced in Section 2, for which we will be able to derive a solution method.

In the subsequent lemmas, we will take the following additional hypothesis:
\begin{itemize}
 \item[$(\mathscr{H})$] $\mu=\pi(\ell,p)$ and $\nu=\pi(\ell,q)$ for some $(\ell,p),(\ell,q)\in\Pi$.
\end{itemize}

\Lem{N1}{Assume $(\mathscr{H})$, let $n\geq2$ and $(f,g),(h,k)\in\C_n(I)$ and define $V:=\big|W_{f,g}^{1,0}\big|^{-p}$. Then $V$ is $(n-1)$-times continuously differentiable. If, for all $x\in I$,
\Eq{2d}{
m_{x;f,g;\mu}''(0)=m_{x;h,k;\nu}''(0)
}
holds, then
\Eq{Phi}{
q\Phi_{h,k}=p\Phi_{f,g}=-\frac{V'}{V}=:\Phi.
}
}

\begin{proof} If $(f,g),(h,k)\in\C_n(I)$, then $W_{f,g}^{1,0}$ is a nowhere zero function which is $(n-1)$-times continuously differentiable, hence so is $V$.

Due to the formula for the second-order derivative in \prp{der} and \lem{0}, we have that
\Eq{*}{
  m_{x;f,g;\mu}''(0)=\mu_2\varphi_{2;f,g}(x)=\mu_2\Phi_{f,g}=2p\ell\Phi_{f,g},\qquad
  m_{x;h,k;\nu}''(0)=\nu_2\varphi_{2;h,k}(x)=\nu_2\Phi_{h,k}=2q\ell\Phi_{h,k}.
}
Thus, the equality \eq{2d} yields the first equality in \eq{Phi}. From the definition of $V$, we get that 
\Eq{*}{
  V'=(-p)\big|W_{f,g}^{1,0}\big|^{-p}\frac{W_{f,g}^{2,0}}{W_{f,g}^{1,0}}
  =(-V)p\Phi_{f,g}=(-V)\Phi.
}
This proves the second equality in \eq{Phi}.
\end{proof}

As an easy consequence of the equality $V'=-V\Phi$, we get the following statement.

\Lem{V}{Under the notations of the previous lemma and appropriate differentiability assumptions, we have
\Eq{*}{
  V'&=-V\Phi,\\
  V''&=V(\Phi^2-\Phi'),\\
  V'''&=V(-\Phi^3+3\Phi'\Phi-\Phi''),\\
  V^{(4)}&=V(\Phi^4-6\Phi'\Phi^2+4\Phi''\Phi+3\Phi'^2-\Phi'''),\\
  V^{(5)}&=V(-\Phi^5+10\Phi'\Phi^3-15\Phi'^2\Phi-10\Phi''\Phi^2+10\Phi''\Phi'+5\Phi'''\Phi-\Phi^{(4)}).
}}

\begin{proof} Using the equality $V'=-V\Phi$, we recursively get
\Eq{*}{
 V''&=(-V\Phi)'=-V'\Phi-V\Phi'
  =V(\Phi^2-\Phi'),\\
 V'''&=(V(\Phi^2-\Phi'))'
  =V'(\Phi^2-\Phi')+V(\Phi^2-\Phi')'
  =V(-\Phi^3+3\Phi'\Phi-\Phi''),\\
 V^{(4)}&=(V(-\Phi^3+3\Phi'\Phi-\Phi''))'
  =V'(-\Phi^3+3\Phi'\Phi-\Phi'')+V(-\Phi^3+3\Phi'\Phi-\Phi'')'\\
  &=V(\Phi^4-6\Phi'\Phi^2+4\Phi''\Phi+3\Phi'^2-\Phi'''),\\
 V^{(5)}&=(V(\Phi^4-6\Phi'\Phi^2+4\Phi''\Phi+3\Phi'^2-\Phi'''))'\\
  &=V'(\Phi^4-6\Phi'\Phi^2+4\Phi''\Phi+3\Phi'^2-\Phi''')+V(\Phi^4-6\Phi'\Phi^2+4\Phi''\Phi+3\Phi'^2-\Phi''')'\\
  &=V(-\Phi^5+10\Phi'\Phi^3-15\Phi'^2\Phi-10\Phi''\Phi^2+10\Phi''\Phi'+5\Phi'''\Phi-\Phi^{(4)}).
}
\end{proof}

\Lem{N2}{
Assume $(\mathscr{H})$, let $n\geq 4$ and $(f,g),(h,k)\in\C_n(I)$ and define $V:=\big|W_{f,g}^{1,0}\big|^{-p}$. If, for all $x\in I$ and $i\in\{1,2\}$,
\Eq{2+4d}{
m_{x;f,g;\mu}^{(2i)}(0)=m_{x;h,k;\nu}^{(2i)}(0)
}
is satisfied, then \eq{Phi} holds and there exist a constant $c\in\R$ and an $(n-3)$-times continuously differentiable function $B:I\to\R$ such that 
\Eq{R}{
\Psi_{f,g}=\frac{B+c}{6p^{\langle 2\rangle}V}-\frac{(p-2)(\Phi'-\Phi^2)}{6p^2}
    \qquad\mbox{and}\qquad
\Psi_{h,k}=\frac{B-c}{6q^{\langle 2\rangle}V}-\frac{(q-2)(\Phi'-\Phi^2)}{6q^2}.
}}

\begin{proof} 
It follows from \lem{N1} that $V$ is $(n-1)$-times, i.e., at least $3$-times continuously differentiable and the case $i=1$ in \eq{2+4d} that \eq{Phi} holds on $I$. Using the recursion \eq{Rec}, we obtain
\Eq{234}{
  \varphi_{2;f,g}&=\Phi_{f,g}=\frac{\Phi}{p}, \qquad
  \psi_{2;f,g}=\Psi_{f,g}, \\
\varphi_{3;f,g}
  &=\varphi_{2;f,g}'+\varphi_{2;f,g}\Phi_{f,g}+\psi_{2;f,g}
  =\frac{\Phi'}{p}+\frac{\Phi^2}{p^2}+\Psi_{f,g},\\
\psi_{3;f,g}
  &=\varphi_{2;f,g}\Psi_{f,g}+\psi_{2;f,g}'
  =-\frac{V'}{pV}\Psi_{f,g}+\Psi_{f,g}',\\
\varphi_{4;f,g}
  &=\varphi_{3;f,g}'+\varphi_{3;f,g}\Phi_{f,g}+\psi_{3;f,g}
  =\frac{\Phi''}{p}+\frac{3\Phi'\Phi}{p^2}+\frac{\Phi^3}{p^3}-\frac{2V'}{pV}\Psi_{f,g}+2\Psi'_{f,g}.
}

These formulae, together with $\mu=\pi(\ell,p)$ and the formula for the fourth-order derivative by \prp{der}, give
\Eq{*}{
m_{x;f,g;\mu}^{(4)}(0)
&=\mu_4\varphi_{4;f,g}(x)-3\mu_2^2(\varphi_{2;f,g}^3+2\psi_{2;f,g}\varphi_{2;f,g})(x)\\
&=12\ell^2\bigg(-\frac{p^{\langle 2\rangle}V'''}{pV}+\frac{3V''V'}{V^2}-\frac{2V'^3}{V^3}+\frac{2p^{\langle 2\rangle}(V\Psi_{f,g})'}{V}\bigg)(x)\\
&=12\ell^2\bigg(\frac{p\Phi''}{1+p}+\frac{3\Phi'\Phi}{1+p}-\frac{\Phi^3}{1+p}+\frac{2p^2}{1+p}\frac{(V\Psi_{f,g})'}{V}\bigg)(x).
}
Similarly, using that $\nu=\pi(\ell,q)$, we get
\Eq{*}{
m_{x;h,k;\nu}^{(4)}(0)
=12\ell^2\bigg(\frac{q\Phi''}{1+q}+\frac{3\Phi'\Phi}{1+q}-\frac{\Phi^3}{1+q}+\frac{2q^2}{1+q}\frac{(V\Psi_{h,k})'}{V}\bigg)(x).
}
Therefore, the case $i=2$ in \eq{2+4d} reduces to
\Eq{*}{
2p^{\langle 2\rangle}(V\Psi_{f,g})'+\frac{p^{\langle 2\rangle}}{p}V(\Phi''-3\Phi'\Phi+\Phi^3)
=2q^{\langle 2\rangle}(V\Psi_{h,k})'+\frac{q^{\langle 2\rangle}}{q}V(\Phi''-3\Phi'\Phi+\Phi^3),
}
or equivalently, using that $V'=-V\Phi$,
\Eq{*}{
2p^{\langle 2\rangle}(V\Psi_{f,g})'+\frac{p^{\langle 2\rangle}}{p}(V(\Phi'-\Phi^2))'
=2q^{\langle 2\rangle}(V\Psi_{h,k})'+\frac{q^{\langle 2\rangle}}{q}(V(\Phi'-\Phi^2))'.
}
Hence, after integration, for some real constant $c$,
\Eq{*}{
2p^{\langle 2\rangle}V\Psi_{f,g}+\frac{p^{\langle 2\rangle}}{p}V(\Phi'-\Phi^2)-\frac{c}{3}
=2q^{\langle 2\rangle}V\Psi_{h,k}+\frac{q^{\langle 2\rangle}}{q}V(\Phi'-\Phi^2)+\frac{c}{3}.
}
Denote by $C$ the function standing on the left hand side of this equality. The regularity assumptions imply that $C$ is $(n-3)$-times continuously differentiable. 
Now defining $B$ as $B:=3C-2V(\Phi'-\Phi^2)$, we can see that $B$ is also $(n-3)$-times continuously differentiable and the above equality yields \eq{R}.
\end{proof}

\Lem{N3}{
Assume $(\mathscr{H})$, let $n\geq 6$ and $(f,g),(h,k)\in\C_n(I)$ and define $V:=\big|W_{f,g}^{1,0}\big|^{-p}$. If, for all $x\in I$ and $i\in\{1,2,3\}$, the equality \eq{2+4d} is satisfied, then the equality \eq{Phi} holds and the identities \eq{R} are valid for some real constant $c$ and an $(n-3)$-times continuously differentiable function $B:I\to\R$. In addition, we get
\Eq{BV}{
(p-q)\bigg(\frac{5pq-p-q-4}{(1+p)(1+q)}V^{(5)}V+(2B''V+B'V')'+B'B\bigg)=c(4pq+3p+3q+2)B'.
}
In particular, if $p=q$, then $cB$ is a constant and if $p\neq q$ and either $5pq=p+q+4$ or $V$ is an at most 4-degree polynomial, then there exist real constants $d,e$ such that
\Eq{BV+}{
(B')^2V
=-\frac{B^3}{6}+\frac{c(4pq+3p+3q+2)}{2(p-q)}B^2+dB+e.
}
}
\begin{proof}
The equalities \eq{Phi} and \eq{R} are consequences of \lem{N1} and \lem{N2}, respectively. Using the recursive formulas \eq{Rec}, \eq{DRec}, the formulas from \eq{234}, 
and the first identity in \eq{R}, we get
\Eq{46}{
\psi_{4;f,g}
  &=\varphi_{3;f,g}\Psi_{f,g}+\psi_{3;f,g}'\\
  &=\bigg(\frac{2\Phi'}{p}+\frac{\Phi^2}{p^2}\bigg)\Psi_{f,g}+\Psi_{f,g}^2+\frac{\Phi}{p}\Psi_{f,g}'+\Psi_{f,g}'',\\
\varphi_{6;f,g}
  &=\varphi_{4;f,g}''+2\varphi_{4;f,g}'\Phi_{f,g}+\varphi_{4;f,g}\varphi_{3;f,g}+2\psi_{4;f,g}'+\psi_{4;f,g}\Phi_{f,g}\\
  &=\frac{\Phi^{(4)}}{p}+\frac{5\Phi'''\Phi}{p^2}+\frac{10\Phi''\Phi'}{p^2}+\frac{10\Phi''\Phi^2}{p^3}+\frac{15\Phi'^2\Phi}{p^3}+\frac{10\Phi'\Phi^3}{p^4}+\frac{\Phi^5}{p^5}
  +\frac{3\Phi}{p}\Psi_{f,g}^2\\
  &\quad+\bigg(\frac{7\Phi''}{p}+\frac{15\Phi'\Phi}{p^2}+\frac{4\Phi^3}{p^3}\bigg)\Psi_{f,g}+\bigg(\frac{12\Phi'}{p}+\frac{9\Phi^2}{p^2}\bigg)\Psi_{f,g}'+6\Psi_{f,g}'\Psi_{f,g}+\frac{9\Phi}{p}\Psi_{f,g}''+4\Psi_{f,g}'''.
}
Combining these identities with $\mu=\pi(\ell,p)$ and the formula for the sixth-order derivative by \prp{der} yields
\Eq{6}{
m_{x;f,g;\mu}^{(6)}(0)
&=\mu_6\varphi_{6;f,g}(x)
-15\mu_4\mu_2(\varphi_{4;f,g}(\varphi_{2;f,g}^2+\psi_{2;f,g})+\varphi_{2;f,g}\psi_{4;f,g})(x)\\
&\quad-15\mu_2^3\varphi_{2;f,g}(\varphi_{3;f,g}\varphi_{2;f,g}^2
-3(\varphi_{2;f,g}^2+\psi_{2;f,g})(\varphi_{2;f,g}^2+2\psi_{2;f,g}))(x)\\
&=120\ell^3p^{\langle3\rangle}\Big(\frac{\Phi^{(4)}}{p}+\frac{5\Phi'''\Phi}{p^2}+\frac{10\Phi''\Phi'}{p^2}+\frac{(7-6p)\Phi''\Phi^2}{p^3}+\frac{15\Phi'^2\Phi}{p^3}
-\frac{(2p+21)\Phi'\Phi^3}{p^3}+\frac{4\Phi^5}{p^3}\\
&\quad+12p\Phi\Psi_{f,g}^2+\frac{2\Psi_{f,g}}{p}((2-3p)\Phi''-15\Phi'\Phi+8\Phi^3)+\frac{6\Psi_{f,g}'}{p}(2\Phi'-3\Phi^2)-12p\Psi_{f,g}'\Psi_{f,g}\\
&\quad+\frac{6(1-p)\Psi_{f,g}''}{p}\Phi+4\Psi_{f,g}'''\Big)(x).
}
Using the first identity in \eq{R} and $V'=-V\Phi$, we recursively get
\Eq{123}{
\Psi'_{f,g}
&=-\frac{p-2}{6p^2}(\Phi''-2\Phi'\Phi)+\frac{B'+(B+c)\Phi}{6p^{\langle2\rangle}V},\\
\Psi''_{f,g}
&=-\frac{p-2}{6p^2}(\Phi'''-2\Phi''\Phi-2\Phi'^2)+\frac{B''+2B'\Phi+(B+c)(\Phi'+\Phi^2)}{6p^{\langle2\rangle}V},\\
\Psi'''_{f,g}
&=-\frac{p-2}{6p^2}(\Phi^{(4)}-2\Phi'''\Phi-6\Phi''\Phi')+\frac{B'''+3B''\Phi+3B'(\Phi'+\Phi^2)+(B+c)(\Phi''+3\Phi'\Phi+\Phi^3)}{6p^{\langle2\rangle}V}.
}
Now, substituting these identities into \eq{6}, we arrive at
\Eq{*}{
m_{x;f,g;\mu}^{(6)}(0)
&=\frac{40\ell^3}{(1+p)(1+2p)}\bigg(p(p+4)(\Phi^{(4)}+\Phi^5)+(7p^2-2p+6)(\Phi^{(3)}\Phi+2\Phi''\Phi')\\
&\quad\quad-(8p^2-13p+9)\Phi''\Phi^2-3(p^2-11p+3)\Phi'^2\Phi-(4p^2+31p-3)\Phi'\Phi^3\\
&\quad\quad+\frac{(p+1)}{V}\Big(2pB'''+(7p+4)B'\Phi'\Big)+\frac{(p+1)^2}{V}\Big(3B''\Phi-B'\Phi^2-\frac{B'(B+c)}{V}\Big)\bigg)(x).
}
Similar argument applies to the case $\nu=\pi(\ell,q)$, we have
\Eq{*}{
m_{x;h,k;\nu}^{(6)}(0)
&=\frac{40\ell^3}{(1+q)(1+2q)}\bigg(q(q+4)(\Phi^{(4)}+\Phi^5)+(7q^2-2q+6)(\Phi^{(3)}\Phi+2\Phi''\Phi')\\
&\quad\quad-(8q^2-13q+9)\Phi''\Phi^2-3(q^2-11q+3)\Phi'^2\Phi-(4q^2+31q-3)\Phi'\Phi^3\\
&\quad\quad+\frac{(q+1)}{V}\Big(2qB'''+(7q+4)B'\Phi'\Big)+\frac{(q+1)^2}{V}\Big(3B''\Phi-B'\Phi^2-\frac{B'(B-c)}{V}\Big)\bigg)(x).
}
Now the case $i=3$ in \eq{2+4d} simplifies to 
\Eq{*}{
(p-q)\bigg(\frac{(5pq-p-q-4)}{(p+1)(q+1)}&V^2\Big(-(\Phi^{(4)}+\Phi^5)+5(\Phi^{(3)}\Phi+2\Phi''\Phi')-10\Phi''\Phi^2-15\Phi'^2\Phi+10\Phi'\Phi^3\Big)\\
&+2B'''-3B''\Phi+B'(\Phi^2-\Phi')+B'B\bigg)=c(4pq+3p+3q+2)B'.
}
Using the first, second and last identities from \lem{V}, we can easily see that 
this equality is exactly equivalent to \eq{BV}.

It is easily seen that the case $p=q$ implies that $cB$ is a constant. On the other hand, if either $5pq=p+q+4$ or $V$ is an at most 4-degree polynomial, then identity \eq{BV} reduces to
\Eq{*}{
(2B''V+B'V')'+B'B=\frac{c(4pq+3p+3q+2)B'}{(p-q)}.
}
Thus, after integration for some real constant $d$, we get the following first-order inhomogeneous linear differential equation for the function $V$: 
\Eq{B''}{
2B''V+B'V'+\frac{B^2}{2}=\frac{c(4pq+3p+3q+2)B}{(p-q)}+d.
}
Multiplying this equality by $B'$ side by side, the equation so obtained is again integrable, therefore there exists a constant $e$ such that formula \eq{BV+} holds.
\end{proof}

\Lem{B''}{
Under the same assumptions of the previous lemma provided that the function $B:I\to\R$ is 5-times continuously differentiable and if $p\neq q$ and either $5pq=p+q+4$ or 
$V$ is an at most 4-degree polynomial. Then we have
\Eq{*}{
  B''V&=\frac{B'V\Phi}{2}+R\circ B,\\
  B'''V&=\frac{B'V}{4}(2\Phi'+\Phi^2)+\frac{3\Phi}{2}(R\circ B)+B'(R'\circ B),\\
  B^{(4)}V
  &=\frac{B'V}{8}\big(4\Phi''+6\Phi'\Phi+\Phi^3\big)
  +\Big(\frac{7\Phi^2}{4}+2\Phi'\Big)(R\circ B)
  +3B'\Phi(R'\circ B)+\frac1{V}((R'R)\circ B)\\
  &\quad+B'^2(R''\circ B),\\
  B^{(5)}V
  &=\frac{B'V}{16}\big(8\Phi'''+16\Phi''\Phi+12\Phi'^2+12\Phi'\Phi^2+\Phi^4\big)+\frac{1}{8}(20\Phi''+50\Phi'\Phi+15\Phi^3)(R\circ B)\\
  &\quad+\frac{B'}{4}\big(25\Phi^2+20\Phi'\big)(R'\circ B)+\frac{5\Phi}{V}((R'R)\circ B)+\frac{B'}{V}(R'^2\circ B)\\
  &\quad+\frac{3B'}{V}(R''R\circ B)+5B'^2\Phi(R''\circ B),
}
where R of at most second degree polynomial which takes the form
\Eq{*}{
R(u):=-\frac{u^2}{4}+\frac{c(4pq+3p+3q+2)u}{2(p-q)}+\frac{d}2.
}
}
\begin{proof}
 Using \lem{N3}, we get that identity \eq{B''} is valid. Thus, using $V'=-V\Phi$, we obtain
\Eq{*}{
  B''V&=\frac{B'V\Phi}{2}-\frac{B^2}{4}-\frac{c(4pq+3p+3q+2)B}{2(p-q)}+\frac{d}2
  =\frac{B'V\Phi}{2}+R\circ B,\\
  B'''V&=-B''V'+\frac{B''V\Phi}{2}+\frac{B'V'\Phi}{2}+\frac{B'V\Phi'}{2}+(R'\circ B)B'\\
  &=\frac{B'V}{4}(2\Phi'+\Phi^2)+\frac{3\Phi}{2}(R\circ B)+B'(R'\circ B),\\
  B^{(4)}V
  &=-B'''V'+\frac{B''V}{4}(2\Phi'+\Phi^2)+\frac{B'V'}{4}(2\Phi'+\Phi^2)+\frac{B'V}{4}(2\Phi''+2\Phi\Phi')\\
  &\quad+\frac{3\Phi'}{2}(R\circ B)+\frac{3\Phi}{2}(R'\circ B)B'+
  (R''\circ B)B'^2+(R'\circ B)B''\\
  &=\Big(\frac{B'V}{4}(2\Phi'+\Phi^2)+\frac{3\Phi}{2}(R\circ B)+(R'\circ B)B'\Big)\Phi \\
  &\quad+\Big(\frac{B'V\Phi}{8}+\frac{R\circ B}{4}\Big)(2\Phi'+\Phi^2)-\frac{B'V\Phi}{4}(2\Phi'+\Phi^2)+\frac{B'V}{4}(2\Phi''+2\Phi\Phi')\\
  &\quad+\frac{3\Phi'}{2}(R\circ B)+\frac{3\Phi}{2}(R'\circ B)B'+
  (R''\circ B)B'^2+\frac{R'\circ B}{V}\Big(\frac{B'V\Phi}{2}+R\circ B\Big)\\
  &=\frac{B'V}{8}\big(4\Phi''+6\Phi'\Phi+\Phi^3\big)
  +\Big(\frac{7\Phi^2}{4}+2\Phi'\Big)(R\circ B)
  +3B'\Phi(R'\circ B)+\frac1{V}((R'R)\circ B)\\
  &\quad+B'^2(R''\circ B),
}
and
\Eq{*}{
  B^{(5)}V
  &=-B^{(4)}V'+\frac{B''V+B'V'}{8}\big(4\Phi''+6\Phi'\Phi+\Phi^3\big)+\frac{B'V}{8}\big(4\Phi'''+6\Phi''\Phi+6\Phi'^2+3\Phi'\Phi^2\big)\\
  &\quad+\Big(\frac{7\Phi'\Phi}{2}+2\Phi''\Big)(R\circ B)
  +\Big(\frac{7\Phi^2}{4}+2\Phi'\Big)(R'\circ B)B'
  +3(B''\Phi+B'\Phi')(R'\circ B)\\&\quad+3B'\Phi(R''\circ B)B'
  -\frac{V'}{V^2}((R'R)\circ B)+\frac1{V}((R''R+R'^2)\circ B)B'+2B''B'(R''\circ B)\\
  &=\Big(\frac{B'V}{8}\big(4\Phi''+6\Phi'\Phi+\Phi^3\big)
  +\Big(\frac{7\Phi^2}{4}+2\Phi'\Big)(R\circ B)
  +3B'\Phi(R'\circ B)+\frac1{V}((R'R)\circ B)\\
  &\quad+B'^2(R''\circ B)\Big)\Phi+\frac{2R\circ B-B'V\Phi}{16}\big(4\Phi''+6\Phi'\Phi+\Phi^3\big)+\frac{B'V}{8}\big(4\Phi'''+6\Phi''\Phi+6\Phi'^2+3\Phi'\Phi^2\big)\\
  &\quad+\Big(\frac{7\Phi'\Phi}{2}+2\Phi''\Big)(R\circ B)
  +\Big(\frac{7\Phi^2}{4}+2\Phi'\Big)(R'\circ B)B'+3\Big(\frac{\Phi}{V}\Big(\frac{B'V\Phi}{2}+R\circ B\Big)+B'\Phi'\Big)(R'\circ B)\\
  &\quad+3B'\Phi(R''\circ B)B'
  +\frac{\Phi}{V}((R'R)\circ B)+\frac{B'}{V}((R''R+R'^2)\circ B)+\frac{2B'}{V}\Big(\frac{B'V\Phi}{2}+R\circ B\Big)(R''\circ B)\\
  &=\frac{B'V}{16}\big(8\Phi'''+16\Phi''\Phi+12\Phi'^2+12\Phi'\Phi^2+\Phi^4\big)+\frac{1}{8}(20\Phi''+50\Phi'\Phi+15\Phi^3)(R\circ B)\\
  &\quad+\frac{B'}{4}\big(25\Phi^2+20\Phi'\big)(R'\circ B)+\frac{5\Phi}{V}((R'R)\circ B)+\frac{B'}{V}(R'^2\circ B)+\frac{3B'}{V}(R''R\circ B)\\
  &\quad+5B'^2\Phi(R''\circ B).
}
\end{proof}

\Lem{N4}{
Assume $(\mathscr{H})$ with $(p,q)=(2,\frac{2}{3})$, let $(f,g),(h,k)\in\C_8(I)$. If, for all $x\in I$ and $i\in\{1,2,3,4\}$, the equality \eq{2+4d} is satisfied, then $\Phi_{h,k}=3\Phi_{f,g}$ holds and there exist real constants $a,b$ such that
\Eq{ab}{
 \Psi_{f,g}=a\Big(W_{f,g}^{1,0}\Big)^2\qquad\mbox{and}\qquad \Psi_{h,k}=b\Big(W_{h,k}^{1,0}\Big)^{\frac{2}{3}}+\frac13\Phi_{h,k}'-\frac{2}{9}\Phi_{h,k}^2.
}}

\begin{proof}
The equality $\Phi_{h,k}=3\Phi_{f,g}$ follows from \lem{N1} and we also have \eq{Phi} with $V:=\big(W_{f,g}^{1,0}\big)^{-2}$. Then $V$ is seven times differentiable and hence $\Phi$ is six times differentiable. The validity of \eq{2+4d} for $i=2$ and \lem{N2} imply the existence of a constant $c$ and a five times differentiable function $B:I\to\R$ such that the equalities in \eq{R} hold.
The parameters $p=2$ and $q=\frac{2}{3}$ also satisfy the condition $5pq=p+q+4$, hence validity of \eq{2+4d} for $i=3$ yields that we have the conclusions of \lem{B''} with  
\Eq{R+}{
  R(u)=-\frac{u^2}{4}+\frac{23cu}{4}+\frac{d}2.
}

Using the two-step recursion \eq{DRec} for $i\in\{4,6\}$, we arrive at
\Eq{*}{
\psi_{6;f,g}
  &=2\varphi_{4;f,g}'\Psi_{f,g}+\varphi_{4;f,g}\psi_{3;f,g}+\psi_{4;f,g}''+\psi_{4;f,g}\Psi_{f,g}\\
  &=\Psi_{f,g}^{(4)}+\Psi_{f,g}'''\frac{\Phi}{p}+7\Psi_{f,g}''\Psi_{f,g}+\Psi_{f,g}''\bigg(\frac{\Phi^2}{p^2}+\frac{4\Phi'}{p}\bigg)
  +\Psi_{f,g}'\bigg(\frac{6\Phi''}{p}+\frac{7\Phi'\Phi}{p^2}+\frac{\Phi^3}{p^3}\bigg)+4\Psi_{f,g}'^2\\
  &\quad+\Psi_{f,g}'\Psi_{f,g}\frac{9\Phi}{p}+\Psi_{f,g}\bigg(\frac{4\Phi'''}{p}+\frac{9\Phi''\Phi}{p^2}+\frac{8\Phi'^2}{p^2}+\frac{9\Phi'\Phi^2}{p^3}
     +\frac{\Phi^4}{p^4}\bigg)+\Psi_{f,g}^2\bigg(\frac{6\Phi'}{p}+\frac{3\Phi^2}{p^2}\bigg)+\Psi_{f,g}^3,
}
and
\Eq{*}{
\varphi_{8;f,g}
  &=\varphi_{6;f,g}''+2\varphi_{6;f,g}'\Phi_{f,g}+\varphi_{6;f,g}\varphi_{3;f,g}+2\psi_{6;f,g}'+\psi_{6;f,g}\Phi_{f,g}\\
  &=\frac{\Phi^{(6)}}{p}+\frac{7\Phi^{(5)}\Phi}{p^2}+\Phi^{(4)}\bigg(\frac{12\Phi'}{p^2}+\frac{21\Phi^2}{p^3}\bigg)
     +35\Phi'''\bigg(\frac{\Phi''}{p^2}+\frac{3\Phi'\Phi}{p^3}+\frac{\Phi^3}{p^4}\bigg)+\frac{70\Phi''^2\Phi}{p^3}\\
     &\quad+35\Phi''\bigg(\frac{3\Phi'^2}{p^3}+\frac{6\Phi'\Phi^2}{p^4}+\frac{\Phi^4}{p^5}\bigg)+\frac{105\Phi'^3\Phi}{p^4}+\frac{105\Phi'^2\Phi^3}{p^5}+\frac{21\Phi'\Phi^5}{p^6}
     +\frac{\Phi^7}{p^7}+\frac{4\Phi\Psi_{f,g}^3}{p}\\
     &\quad+\bigg(\frac{22\Phi''}{p}+\frac{42\Phi'\Phi}{p^2}+\frac{10\Phi^3}{p^3}\bigg)\Psi_{f,g}^2+\bigg(\frac{16\Phi^{(4)}}{p}+\frac{56\Phi'''\Phi}{p^2}+\frac{86\Phi''\Phi^2}{p^3}
     +\frac{112\Phi''\Phi'}{p^2}+\frac{128\Phi'^2\Phi}{p^3}\\
     &\quad+\frac{70\Phi'\Phi^3}{p^4}+\frac{6\Phi^5}{p^5}\bigg)\Psi_{f,g}+\bigg(\frac{46\Phi'''}{p}+\frac{124\Phi''\Phi}{p^2}+\frac{90\Phi'^2}{p^2}+\frac{142\Phi'\Phi^2}{p^3}
     +\frac{20\Phi^4}{p^4}\bigg)\Psi_{f,g}'\\
     &\quad+\frac{40\Phi}{p}\Psi_{f,g}'^2+\bigg(\frac{72\Phi'}{p}+\frac{48\Phi^2}{p^2}\bigg)\Psi_{f,g}'\Psi_{f,g}+\bigg(\frac{60\Phi''}{p}+\frac{124\Phi'\Phi}{p^2}+\frac{34\Phi^3}{p^3}\bigg)\Psi_{f,g}''
     +\frac{52\Phi}{p}\Psi_{f,g}''\Psi_{f,g}\\
     &\quad+\bigg(\frac{44\Phi'}{p}+\frac{34\Phi^2}{p^2}\bigg)\Psi_{f,g}'''+\frac{20\Phi}{p}\Psi_{f,g}^{(4)}+12\Psi_{f,g}'\Psi_{f,g}^2+48\Psi_{f,g}''\Psi_{f,g}'
     +24\Psi_{f,g}'''\Psi_{f,g}+6\Psi_{f,g}^{(5)}.
}
These identities, together with the formulas from \eq{234} and \eq{46}, $\mu=\pi(\ell,2)$, and the formula for the eighth-order derivative by \prp{der}, imply that 
\Eq{8}{
 m_{x;f,g;\mu}^{(8)}(0)&=256\ell^4\Big(\varphi_{8;f,g}(x)-28(\varphi_{6;f,g}(\varphi_{2;f,g}^2+\psi_{2;f,g})
  +\varphi_{2;f,g}\psi_{6;f,g})(x)-35(\varphi_{4;f,g}^2\varphi_{2;f,g}\\
  &\quad+2\varphi_{4;f,g}\psi_{4;f,g})(x)+210(\varphi_{4;f,g}(3\varphi_{2;f,g}^4+\varphi_{2;f,g}^2(7\psi_{2;f,g}-\varphi_{3;f,g})+2\psi_{2;f,g}^2)\\
  &\quad+2\varphi_{2;f,g}\psi_{4;f,g}(\varphi_{2;f,g}^2+2\psi_{2;f,g}))(x)-105(\varphi_{4;f,g}\varphi_{2;f,g}^4+15\varphi_{2;f,g}^7\\
  &\quad+2\varphi_{2;f,g}^3(5\varphi_{2;f,g}^2+6\psi_{2;f,g})(6\psi_{2;f,g}-\varphi_{3;f,g})+4\varphi_{2;f,g}(6\psi_{2;f,g}^3-\varphi_{2;f,g}^3\psi_{3;f,g}))(x)\Big)\\
  &=16\ell^4\Big(8\Phi^{(6)}+28\Phi^{(5)}\Phi+14\Phi^{(4)}(6\Phi'-\Phi^2)+35\Phi'''(4\Phi''+6\Phi'\Phi-3\Phi^3)+70\Phi''^2\Phi\\
  &\quad+70\Phi''\Phi'(3\Phi'-7\Phi^2)+105\Phi'^3\Phi-630\Phi'^2\Phi^3+329\Phi'\Phi^5-34\Phi^7-8704\Phi\Psi_{f,g}^3\\
  &\quad+64(22\Phi''+119\Phi'\Phi-85\Phi^3)\Psi_{f,g}^2-16(6\Phi^{(4)}+49\Phi'''\Phi+77\Phi''\Phi'-51\Phi''\Phi^2+117\Phi'^2\Phi\\
  &\quad-238\Phi'\Phi^3+51\Phi^5)\Psi_{f,g}+8(46\Phi'''-127\Phi''\Phi+45\Phi'^2+360\Phi'\Phi^2+125\Phi^4)\Psi_{f,g}'\\
  &\quad-4352(\Phi'-2\Phi^2)\Psi_{f,g}'\Psi_{f,g}-2816\Phi\Psi_{f,g}'^2+2432\Phi\Psi_{f,g}''\Psi_{f,g}-8(10\Phi''+99\Phi'\Phi-26\Phi^3)\Psi_{f,g}''\\
  &\quad+8(44\Phi'-53\Phi^2)\Psi_{f,g}'''-64\Phi\Psi_{f,g}^{(4)}+32((272\Psi_{f,g}^2-46\Psi_{f,g}'')\Psi_{f,g}'-44\Psi_{f,g}'''\Psi_{f,g}+3\Psi_{f,g}^{(5)})\Big)(x).
}
Using the identities \eq{123} and $V'=-V\Phi$, we get
\Eq{*}{
\Psi^{(4)}_{f,g}
&=-\frac{p-2}{6p^2}(\Phi^{(5)}-2\Phi^{(4)}\Phi-8\Phi'''\Phi'-6\Phi''^2)+\frac{1}{6p^{\langle2\rangle}V}\Big(B^{(4)}+4B'''\Phi+6B''(\Phi'+\Phi^2)\\
&\quad+3B'(\Phi''+3\Phi'\Phi+\Phi^3)+(B+c)(\Phi'''+4\Phi''\Phi+6\Phi'\Phi^2+3\Phi'^2+\Phi^4)\Big),\\
\Psi^{(5)}_{f,g}
&=-\frac{p-2}{6p^2}(\Phi^{(6)}-2\Phi^{(5)}\Phi-10\Phi^{(4)}\Phi'-20\Phi'''\Phi'')+\frac{1}{6p^{\langle2\rangle}V}\Big(B^{(5)}+4B^{(4)}\Phi+10B'''(\Phi'+\Phi^2)\\
&\quad+9B''(\Phi''+3\Phi'\Phi+\Phi^3)+B'(4\Phi'''+16\Phi''\Phi+12\Phi'^2+21\Phi'\Phi^2+4\Phi^4)\\
&\quad+(B+c)(\Phi^{(4)}+5\Phi'''\Phi+10\Phi''\Phi'+10\Phi''\Phi^2+15\Phi'^2\Phi+10\Phi'\Phi^3+\Phi^5)\Big).
}
Substituting these identities and \eq{123} into \eq{8}, we have
\Eq{*}{
m_{x;f,g;\mu}^{(8)}(0)&=16\ell^4\Bigg(8\Phi^{(6)}+28\Phi^{(5)}\Phi+14\Phi^{(4)}(6\Phi'-\Phi^2)+35\Phi'''(4\Phi''+6\Phi'\Phi-3\Phi^3)+70\Phi''^2\Phi\\
&\quad+70\Phi''\Phi'(3\Phi'-7\Phi^2)+105\Phi'^3\Phi-630\Phi'^2\Phi^3+329\Phi'\Phi^5-34\Phi^7+\frac{B'''(164\Phi'+35\Phi^2)}{V}\\
&\quad+\frac{1}{V}\Big(B''(110\Phi''+345\Phi'\Phi-61\Phi^3)+B'(106\Phi'''+61\Phi''\Phi+357\Phi'^2-321\Phi'\Phi^2+46\Phi^4)\Big)\\
&\quad-\frac{B+c}{V^2}\Big(22B'''+51B''\Phi+B'(157\Phi'+11\Phi^2)-\frac{17B'(B+c)}{V}\Big)+\frac{\Phi}{V}\Big(52B^{(4)}-\frac{90B'^2}{V}\Big)\\
&\quad+\frac{1}{V}\Big(12B^{(5)}-\frac{23B''B'}{V}\Big)\Bigg)(x).
}
Similarly, using that $\nu=\pi(\ell,\frac23)$, we arrive at
\Eq{*}{
m_{x;h,k;\nu}^{(8)}(0)&=16\ell^4\Bigg(8\Phi^{(6)}+28\Phi^{(5)}\Phi+14\Phi^{(4)}(6\Phi'-\Phi^2)+35\Phi'''(4\Phi''+6\Phi'\Phi-3\Phi^3)+70\Phi''^2\Phi\\
&\quad+70\Phi''\Phi'(3\Phi'-7\Phi^2)+105\Phi'^3\Phi-630\Phi'^2\Phi^3+329\Phi'\Phi^5-34\Phi^7+\frac{B'''(205\Phi'+380\Phi^2)}{3V}\\
&\quad+\frac{1}{3V}\Big(5B''(74\Phi''+267\Phi'\Phi-55\Phi^3)+B'(326\Phi'''+187\Phi''\Phi+1139\Phi'^2\\
&\quad\quad-1111\Phi'\Phi^2+162\Phi^4)\Big)-\frac{B-c}{3V^2}\Big(50B'''+225B''\Phi+B'(575\Phi'+25\Phi^2)-\frac{75B'(B-c)}{V}\Big)\\
&\quad+\frac{70\Phi}{3V}\Big(2B^{(4)}-\frac{5B'^2}{V}\Big)+\frac{5}{3V}\Big(4B^{(5)}-\frac{B''B'}{V}\Big)\Bigg)(x).
}
Therefore, the case $i=4$ in \eq{2+4d} reduces to
\Eq{*}{
&16B^{(5)}V+16B^{(4)}V\Phi-4B'''(4B+29c)-64B''B'-\frac{12B'(2B^2-21Bc+2c^2)}{V}+80B'^2\Phi\\
&+B'\Big(B(104\Phi'-8\Phi^2)-c(1046\Phi'+58\Phi^2)\Big)+B''(72B-378c)\Phi+B'''\Big(112\Phi'-100\Phi^2\Big)V\\
&-B''\Big(40\Phi''+300\Phi'\Phi-92\Phi^3\Big)V-B'\Big(8\Phi'''+4\Phi''\Phi+68\Phi'^2-148\Phi'\Phi^2+24\Phi^4\Big)V=0.
}
Using \lem{B''}, this identity simplifies to
\Eq{R'}{
  &\Big(-16B'+6\Phi(2B-23c)\Big)(R\circ B)+B'\Big(12V(4\Phi'+\Phi^2)-(4B+29c)\Big)(R'\circ B)\\
  &+24\Phi((R'R)\circ B)+4B'(R'^2\circ B)+12B'(R''R\circ B)+12B'^2V\Phi\Big(2(R''\circ B)+1\Big)\\
  &-3B'(2B^2-21Bc+2c^2)+3B'V(2B-23c)(4\Phi'+\Phi^2)=0.
}
Now, by \eq{R+}, we get
\Eq{*}{
R\circ B= -\frac{B^2}{4}+\frac{23cB}{4}+\frac{d}{2},\qquad R'\circ B=-\frac{B}{2}+\frac{23c}{4}\qquad\mbox{and}\qquad R''\circ B=-\frac{1}{2}.
}
Using these identities, then \eq{R'} reduces to
\Eq{B'}{
B'&\big(5B^2-190cB-81c^2-22d\big)=0.
}
We show that this equality implies that $B$ is constant on $I$. To the contrary, assume that $B$ is nonconstant on $I$. Then $B'$ is not identically zero, and hence $B'(x_0)\neq0$ for some $x_0\in I$. By the continuity of $B'$, there exists an open subinterval $J$ of $I$ containing $x_0$ such that $B'(x)\neq0$ for all $x\in J$. Then, for all $x\in J$, the equality \eq{B'} yields that
\Eq{*}{
  5B^2(x)-190cB(x)-81c^2-22d=0.
}
That is, for all $x\in J$, the value $B(x)$ is equal to one of the roots of the second degree polynomial $5u^2-190cu-81c^2-22d$. By the continuity of $B$, now it follows that $B$ is constant on $J$, which yields $B'(x_0)=0$, contradicting the choice of $x_0$. Therefore, $B$ has to be a constant function. Then denoting 
$\frac{B+c}{6p^{\langle2\rangle}}$ by $a$ and $\frac{B-c}{6q^{\langle2\rangle}}$ by $b$, the equalities in \eq{R} simplify to \eq{ab}. 
\end{proof}

Now we are ready to present and prove our main result.

\Thm{M}{
Assume $(\mathscr{H})$ with $(p,q)=(2,\frac23)$ and let $(f,g),(h,k)\in\C_8(I)$. Then the following assertions are equivalent:
\begin{enumerate}[(i)]
 \item The means $\M_{f,g;\mu}$ and $\M_{h,k;\nu}$ are equal on $I^2$.
 \item The means $\M_{f,g;\mu}$ and $\M_{h,k;\nu}$ are equal near the diagonal of  $I^2$.
 \item For all $x\in I$ and $i\in\{1,2,3,4\}$, the identities $m_{x;f,g;\mu}^{(2i)}(0)=m_{x;h,k;\nu}^{(2i)}(0)$ hold.
 \item $\Phi_{h,k}=3\Phi_{f,g}$ holds on $I$ and there exist real constants $a$ and $b$ such that the equalities in \eq{ab} hold on $I$.
\item There exist real constants $\alpha,\beta,\gamma,\delta,\varepsilon,\zeta,\eta$ such that
\Eq{fghk}{
\alpha f^2+\beta fg+\gamma g^2=1\qquad\mbox{and}\qquad \delta h^2+\varepsilon hk+\zeta k^2=\Big(W_{h,k}^{1,0}\Big)^{\frac{2}{3}}
}
and $W_{h,k}^{1,0}=\eta\Big(W_{f,g}^{1,0}\Big)^3$.
\item There exist two real polynomials $P$ and $Q$ of at most second degree which are positive on the range of $f/g$ and $h/k$, respectively, and there exist real constants $\eta$ and $\rho$ such
\Eq{gk}{
g=\frac{1}{\sqrt{P}}\circ \frac{f}{g},\quad k=\bigg|\bigg(\frac{h}{k}\bigg)'\bigg|\Bigg(\frac{1}{\sqrt{Q^3}}\circ\frac{h}{k}\Bigg), \quad\mbox{and}\quad 
\bigg(\int\frac{1}{Q}\bigg)\circ\frac{h}{k}
   =\eta^{\frac{1}{3}}\bigg(\int\frac{1}{P}\bigg)\circ\frac{f}{g}+\rho.
}
\item $\M_{f,g;\mu}=\A_\varphi=\M_{h,k;\nu}$ holds on $I^2$ with $\varphi:=\int W_{f,g}^{1,0}$.
\item There exists a strictly monotone function $\varphi:I\to\R$ such that $\M_{f,g;\mu}=\A_\varphi=\M_{h,k;\nu}$ holds on $I^2$.
\item There exist a strictly monotone differentiable function $\varphi:I\to\R$ and real constants $a$ and $b$ such that
 \Eq{*}{
   (f,g)\sim(S_a\circ\varphi, C_a\circ\varphi) \qquad\mbox{and}\qquad
   (h,k)\sim(\varphi'\cdot S_b\circ\varphi,\varphi'\cdot C_b\circ\varphi).
 }
\end{enumerate}
}
\begin{proof}
The implication (i)$\Rightarrow$(ii) is obvious. The implications (ii)$\Rightarrow$(iii) and (iii)$\Rightarrow$(iv) are consequences of \lem{N0} and \lem{N4}, respectively.

Assume that (iv) is valid for some real constants $a,b$. Then, the integration of the identity $\Phi_{h,k}=3\Phi_{f,g}$ implies 
\Eq{W+}{
W_{h,k}^{1,0}=\eta\Big(W_{f,g}^{1,0}\Big)^3
}
for some real constant $\eta$. Now the first identity in \eq{ab} and implication 
(iv)$\Rightarrow$(ii) of \cite[Theorem 10]{PalZak20b} yield that there exist real constants $\alpha,\beta,\gamma$ such that the first equality in \eq{fghk} holds. 
An easy computation, using the second identity in \eq{ab}, implies that the expression
\Eq{*}{
\frac{3W^{3,0}_{h,k}+12W^{2,1}_{h,k}}{\Big(W^{1,0}_{h,k}\Big)^{\frac53}}
 -5\frac{\Big(W^{2,0}_{h,k}\Big)^2}{\Big(W^{1,0}_{h,k}\Big)^{\frac83}}
}
is constant. Applying implication (vii)$\Rightarrow$(iv) of \cite[Theorem 7]{LovPalZak20}, then there exist real constants $\delta,\varepsilon,\zeta$ such that the second equality in \eq{fghk} is valid. 
Thus, we have proved that assertion (v) holds.

To prove implication (v)$\Rightarrow$(vi), assume that (v) is valid for some real constants $\alpha,\beta,\gamma,\delta,\varepsilon,\zeta,\eta$ and define
\Eq{*}{
P(u):=\alpha u^2+\beta u+\gamma\qquad\mbox{and}\qquad Q(u):=\delta u^2+\varepsilon u+\zeta\qquad (u\in\R).
}
It is easily seen that the two identities in \eq{fghk} are equivalent to the following two identities
\Eq{PQ}{
P\circ\frac{f}{g}=\frac{1}{g^2}\qquad\mbox{and}\qquad Q\circ\frac{h}{k}=\frac{\Big(W_{h,k}^{1,0}\Big)^{\frac{2}{3}}}{k^2},
}
respectively. Thus, $P$ and $Q$ are positive real polynomials on the range of the ratio functions $f/g$ and $h/k$, respectively. The first equality in \eq{gk} is a result of the first identity in \eq{PQ}. 
The second equality in \eq{gk} is a consequence of $W^{1,0}_{h,k}=k^2(h/k)'$ and the second identity in \eq{PQ}. Furthermore, using \eq{PQ}, we have
that
\Eq{*}{
  \bigg(\bigg(\int\frac{1}{P}\bigg)\circ\frac{f}{g}\bigg)'
  =\bigg(\frac{1}{P}\circ\frac{f}{g}\bigg)\cdot\bigg(\frac{f}{g}\bigg)'
  =g^2\cdot\frac{W^{1,0}_{f,g}}{g^2}
  =W_{f,g}^{1,0}.
}
This identity, together with \eq{PQ} and \eq{W+}, implies
\Eq{*}{
  \bigg(\bigg(\int\frac{1}{Q}\bigg)\circ\frac{h}{k}\bigg)'
  &=\bigg(\frac{1}{Q}\circ\frac{h}{k}\bigg)\cdot\bigg(\frac{h}{k}\bigg)'
  =\frac{k^2}{\Big(W^{1,0}_{h,k}\Big)^{\frac23}}\cdot\frac{W^{1,0}_{h,k}}{k^2}
  =\Big(W_{h,k}^{1,0}\Big)^{\frac13}=\eta^{\frac{1}{3}}W_{f,g}^{1,0}\\
  &=\eta^{\frac{1}{3}}\bigg(\bigg(\int\frac{1}{P}\bigg)\circ\frac{f}{g}\bigg)'.
}
Integrating both sides, then we arrive at the last equality in \eq{gk} for some real constant $\rho$. Therefore, assertion (vi) is valid. 
The implication (vi)$\Rightarrow$(v) is obvious by reversing all the implications in the previous calculation. Thus (vi) and (v) are shown to be equivalent. 

To prove that (v) implies (vii). Assume that (v) holds for some real constants $\alpha,\beta,\gamma,\delta,\varepsilon,\zeta,\eta$. Then implication (ii)$\Rightarrow$(iii) of \cite[Theorem 10]{PalZak20b} 
and implication (iv)$\Rightarrow$(vi) of \cite[Theorem 7]{LovPalZak20} imply that 
\Eq{*}{
\M_{f,g;\mu}=\A_\varphi \qquad\mbox{and}\qquad \M_{h,k;\nu}=\A_\psi,
}
respectively, hold on $I^2$ with $\varphi=\int W_{f,g}^{1,0}$ and $\psi=\int \Big(W_{h,k}^{1,0}\Big)^{\frac13}$. The identity \eq{W+} gives $\psi=\eta^{\frac13}\varphi$, from which we get that $\A_\varphi=\A_\psi$ is valid on $I^2$. Consequently, (vii) holds. 

The implications (vii)$\Rightarrow$(viii) and (viii)$\Rightarrow$(i) are straightforward. Therefore, all the assertions from (i) to (viii) are equivalent. 
Finally, the equivalence of (viii) and (ix) is a direct consequence of \cite[Corollary 9]{PalZak20b} and \cite{LovPalZak20}, respectively. 
\end{proof}

For the proof of the equivalence of the nine conditions of the theorem, we essentially needed the eight times differentiability of the generating functions $f,g,h,k$. On the other hand, one can observe that some particular implications can be verified under weaker regularity assumptions. The assertion (i) trivially follows from (viii) in the case $(f,g),(h,h)\in\C_0(I)$. It is an open problem whether the reversed implication is also true with this natural regularity assumption.

Finally, we can answer the question formulated in the introduction about two-variable means that are equal to a Bajraktarevi\'c mean and to a Cauchy mean at the same time.

\Cor{M}{
Assume $(\mathscr{H})$ with $(p,q)=(2,\frac23)$. Then the intersection of the 
classes of means
\Eq{IS}{
  \{\M_{f,g;\mu}\mid (f,g)\in\C_8(I)\}\,\cap\, 
  \{\M_{h,k;\nu}\mid (h,k)\in\C_8(I)\}
}
consists of the symmetric two-variable quasiarithmetic means. In other words, if a two-variable mean is simultaneously is a Bajraktarevi\'c mean and a Cauchy mean with eight times differentiable generators, then it has to be a quasiarithmetic mean.}

\begin{proof} Assume that a two-variable mean $M:I^2\to I$ belongs to both of the classes of means. Then there exist $(f,g)\in\C_8(I)$ and $(h,k)\in\C_8(I)$ such that $M=\M_{f,g;\mu}$ and $M=\M_{h,k;\nu}$. Hence, assertion (i) of \thm{M} holds. Then, by this theorem assertion (viii) is also valid, hence $M$ has to be quasiarithmetic. 
\end{proof}

As an application of this corollary, we can deduce the result of Alzer and Ruscheweyh.

\Cor{AR}{Assume that a two-variable mean $M:\R_+^2\to\R_+$ is simultaneously equal to a Gini mean and to a Stolarsky mean. Then $M$ has to be a power mean.}

\begin{proof} If $M$ equals to a Gini mean, then there exist real constants $a,b$ such that $M=\G_{a,b}$, where $\G_{a,b}$ was defined by \eq{Gini}. This shows that $M$ is homogeneous and it also equals a Bajraktarevi\'c mean with generators $f$ and $g$ defined by $f(x):=x^a$ and $g(x):=x^b$ if $a\neq b$ and $f(x):=x^a\log(x)$ if $a=b$. 

If $M$ also equals to a Stolarsky mean, then there exist real constants $c,d$ such that $M=\S_{c,d}$, where $\S_{c,d}$ was defined by \eq{Sto}. Then $M$ equals a Cauchy mean with obviously chosen generators. Thus, $M$ belongs to the intersection \eq{IS} and hence it has to be a quasiarithmetic mean. Being also homogeneous, it follows that it has to be a power mean. 
\end{proof}

%\bibliography{publ,funcequ}
%\bibliographystyle{amsplain}

\providecommand{\bysame}{\leavevmode\hbox to3em{\hrulefill}\thinspace}
\providecommand{\MR}{\relax\ifhmode\unskip\space\fi MR }
\providecommand{\MRhref}[2]{%
  \href{http://www.ams.org/mathscinet-getitem?mr=#1}{#2}
}
\providecommand{\href}[2]{#2}

\end{document}